\documentclass[12pt,leqno]{amsart}
\usepackage{amsmath, amsfonts, amssymb,amsthm ,amscd}
\usepackage{graphicx}
\usepackage{epsf}
\usepackage{epsfig}
\input rlepsf
\usepackage{psfrag}
\usepackage{color}

\numberwithin{equation}{section}
\newtheorem{lemma}{Lemma}[section]
\newtheorem{theorem}[lemma]{Theorem}

\newtheorem{proposition}[lemma]{Proposition}
\newtheorem{definition}[lemma]{Definition}

\newtheorem{remark}[lemma]{Remark}
\let\lutzremark=\remark
\def\remark{\lutzremark\normalfont}

\newtheorem{notation}[lemma]{Notation}
\newtheorem{example}[lemma]{Example}

\renewcommand{\epsilon}{\varepsilon}

\newcommand{\pr}{\operatorname{pr}}

\newcommand{\R}{{\mathbb R}}

{\catcode`"=12 \gdef\hex{"}}

\mathchardef\laplace=\hex0001

\mathchardef\nabla=\hex0272

\makeatletter

\def\@@dalembert#1#2{\setbox0\hbox{$#1\mathrm I$}

  \vrule height\ht0 depth\z@ width.04\ht0

  \rlap{\vrule height\ht0 depth-.96\ht0 width.8\ht0}

  \vrule height.1\ht0 depth\z@ width.8\ht0

  \vrule height\ht0 depth\z@ width.1\ht0 }

\def\dalembert{\mathbin{\mathpalette\@@dalembert{}}\,}

\makeatother

\makeatletter

  \def\restrict#1{\ifmmode\def\next{\mathpalette

  \mathrestrict}\fi\left.\next{#1}\right._}

  \def\mathrestrict#1#2{\setbox0=\hbox{$\m@th#1{#2}$}\relax

  \dimen@=\dp0 \advance\dimen@ by 0.7ex\relax

  #2\,\vrule height\ht0 width 0.4pt depth\dimen@}

\makeatother

\makeatletter

\def\@@varcirc#1#2{\mathbin{\lower#1ex\hbox{\m@th${#2\mathchar\hex0017 }$}}}
\def\varcirc{\mathchoice
  {\@@varcirc{.91}\displaystyle}{\@@varcirc{.91}\textstyle}%
{\@@varcirc{.57}\scriptstyle}{\@@varcirc{.45}\scriptscriptstyle}}
\makeatother

\begin{document}

\title[]{A General Fredholm Theory and Applications}
\author[H.~Hofer]{H.~Hofer$^{\dagger}$}
 \address[H.~Hofer]{Courant Institute, 251 Mercer Street, New York, NY 10012, USA}

\thanks{$\dagger$ Research partially supported by  NSF grants DMS-$0102298$ and
DMS-$0505968$}

\maketitle

%\tableofcontents

%\section{Introduction}
The theory described here results from an attempt to find a general
abstract framework in which various theories, like Gromov-Witten
Theory (GW), Floer Theory (FT), Contact Homology (CH) and more
generally Symplectic Field Theory (SFT) can be understood from a
general point of view. Let us describe the general landscape in a
somewhat oversimplified form. The common feature (with the exception
of GW which has less structure) is the fact that we have infinitely
many different Fredholm problems defined on spaces with boundary
with corners, where the boundary strata can be explained in terms of
products (or more generally fibered products) of other problems (on
the list). In oversimplified form, the solution sets are zeros of a
section $f$ of some bundle $\tau:Y\rightarrow X$, where the space
has a boundary $\partial X$, and where moreover there exists a
recipe (or even many recipes) to construct from two given
solutions\footnote{In general one requires them to satisfy a
compatibility condition.} $x'$ and $x''$ of $f=0$ a new solution,
say the product, $x=x'\circ x''$. The recipes for constructing new
solutions are defined even for non-solutions and $\partial X$ is
precisely the space of points which are products. Hence we have
$$
\partial X =X\circ X.
$$
Moreover, if we denote the restriction of $f$ to $\partial X$ by
$\partial f$ and define $f\circ f$ on $\partial X$ as the set-valued
section
$$
f\circ f(x)=\{f(x')\circ f(x'')\ |\ x=x'\circ x''\}
$$
then we say that $f$ is compatible with the recipe $\circ$ provided
$$
\partial f =f\circ f.
$$
Assuming $f$ to be compatible the $\circ$-structure generates a
certain amount of algebra which can be captured on a rather
rudimentary level. Then the more sophisticated algebra we see in the
description of SFT can be viewed as obtained by some kind of
representation theory of the underlying "primitive" data. One can
develop a general theory which covers a variety of problems. The
only difference in application between seemingly different problems
is that the underlying structure of the family of Fredholm problems,
i.e. their interaction, is different.

The starting point for the investigation is the Symplectic Field
Theory (SFT) as initiated by Eliashberg, Givental and the author in
\cite{EGH}. Wysocki, Zehnder and the author developed a powerful
nonlinear Fredholm theory with operations (FTO) which can be used to
describe SFT, \cite{HWZ1,HWZ2,HWZ3}. The Fredholm theory takes place
in a new kind of spaces called polyfolds. These spaces are needed
since all phenomena of interest are coming from analytically
difficult phenomena like bubbling-off, stretching the neck, breaking
of trajectories and blowing-up\footnote{Perhaps folklore-wise known
as the "Analytical Chamber of Horror".}. In \cite{EH} this theory
will be used to develop SFT  and many ideas around it  in full
generality (and absolutely rigorously). In fact, the polyfold
language allows to completely remove the analysis before carrying
out topological and algebraic considerations contrasting the current
situation where many arguments are relying on a combination of
arguments across the board.

Let us discuss for a moment the issues which we have to address in
developing a general framework allowing us to describe seemingly
different problems. Gromov-Witten theory, Floer-Theory,
Contact-Homology, or more generally Symplectic Field Theory are
theories build on the study of certain compactified moduli spaces,
or even infinite families of such spaces. These moduli spaces are
measured and the data is encoded in convenient ways, quite often as
a so-called generating function. Common features include:
\begin{itemize}
\item[1)] The moduli spaces are solutions of elliptic PDE's quite
often exhibiting compactness problems, at least as seen from a
more classical analytical viewpoint.

\item[2)] Very often these moduli spaces, when they are not
compact admit nontrivial compactifications usually based on
surviving analytic phenomena carrying names like "Bubbling-off",
"Stretching the Neck", "Blow-up", "Breaking of Trajectories" hinting
to borderline analytic behavior.

\item[3)] In problems like Floer-Theory, Contact-Homology or
Symplectic Field Theory precisely the algebraic structures of
interest are those created by  the "violent analytic behavior" and
its "taming" by finding a workable compactification. In fact the
algebra is created by the fact that many different moduli spaces
interact with each other in a complicated way.
\end{itemize}

We begin with the shortcomings of classical Fredholm theory. The
classical Fredholm theory can be viewed as the study of Fredholm
sections of some Banach bundle $Y\rightarrow X$. For definiteness we
assume that $Y$ is a Banach space bundle over the Banach manifold
$X$. Let us denote the fiber over $x\in X$ by $Y_x$. If $f(x)=0$ we
can build the linearisation $f'(x):T_xX\rightarrow Y_x$ and if
$f'(x)$ is surjective we have a solution manifold near $x$ in fact
inheriting its manifold structure as  a submanifold of the (big)
ambient space. From a practical point of view the bubbling-off
phenomena usual cannot be described within this classical framework.
The key question is therefore if there is a generalized Fredholm
theory in which interesting problems of the type described above can
be handled. Keeping this in mind it is worthwhile to have a critical
look at the classical case. We may  raise the following question in
the classical context. Is it not "unnecessary luxury" that the
ambient space has a lot of "hard structure" whereas we only seem to
use little of it in order to obtain a smooth structure on the
solution set $f^{-1}(0)$ (assuming transversality)? This question is
very much justified, since in many cases, once the solution spaces
are constructed, the ambient spaces are discarded and considered
irrelevant.  The hope is, of course, that analyzing the situation,
we might be able to see what is the bare minimum of structure needed
for a suitable generalization. More precisely we have to address the
following question:\\
{ What (perhaps new) structures do we need on
the ambient space and bundle (with a preferred section called $0$)
to talk about transversality and an abstract perturbation theory for
a section $f$ so that at points of transversality the solution set
$f^{-1}(0)$ carries in a natural way the structure of a smooth
orbifold with boundary with corners? In addition we require the
theory to be so general that in applications the compactified moduli
spaces in Gromov-Witten theory, Floer theory or
SFT would be the solution sets of the generalized Fredholm operators.}\\
Analyzing the before-mentioned theories it becomes immediately clear
that one has to address a certain number of very serious issues. For
example, if one of the "violent analytical phenomena" occurs, any
natural candidate for an ambient space seems to have  locally
varying dimensions. With other words, in particular, the spaces are
locally not isomorphic to open sets in Banach or Frechet spaces.
Hence, if we still think we should devise  a manifold-type theory,
the local models cannot be open sets in some Banach or Frechet
space. They need to be more general. If we still want to talk about
a linearisation of a problem, which, as every analyst knows, has its
undeniable benefits, we should look for some class of local models
which in some way admit tangent spaces. Moreover, there are some
other unpleasant phenomena to deal with. For example in some
constructions we have divide out  by families of diffeomorphisms
acting on the domain of maps. Analysts  know that such actions (for
example in any Banach space set-up) will always be only continuous,
but never smooth. In addition, the applications like SFT, will
require the theory to have certain features of a theory of
(infinite-dimensional) orbifolds with boundaries with corners.
\begin{figure}[htbp]
\mbox{}\\[2ex]
\centerline{\relabelbox \epsfxsize 5truein \epsfbox{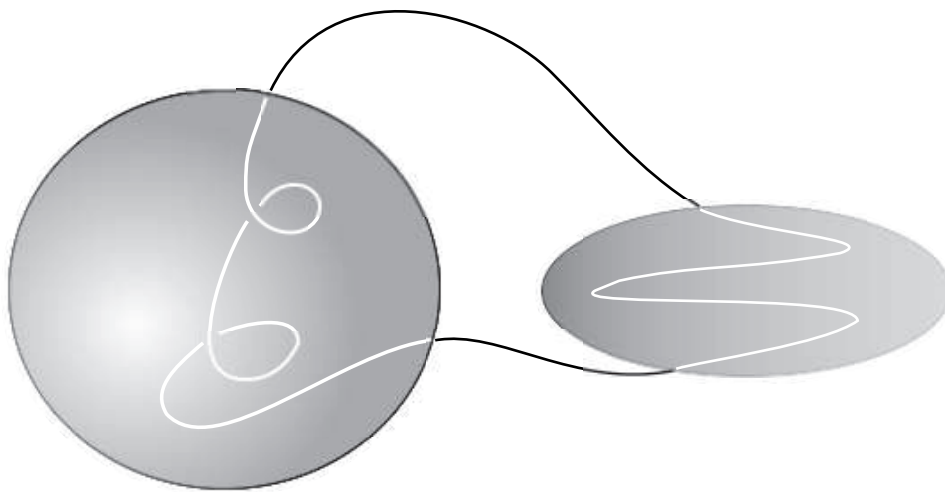}
\endrelabelbox}
\caption{This figure shows a finite-dimensional $M$-polyfold, say $X$, homeomorphic to
the space obtained from an open three-ball and an open two-ball
connected by two curves having a one-dimensional $S^1$-like
submanifold. This submanifold could arise as the zero of a
transversal section of a strong (infinite-dimensional) M-polyfold bundle $Y$ over $X$,
which has varying dimensions. Namely over the three-ball it is
two-dimensional, over the two-disk one-dimensional and otherwise
trivial. The polyfold theory guarantees natural smooth structures on
such solution sets.}\label{Fig1.0}
\mbox{}\\[1ex]
\end{figure}
Summarizing there is a whole basket of issues which call for a more
general theory. If we have a look at our list of requirements it
seems  that the problem of finding an adequate theory is
"over-determined". Surprisingly, however, there is such a general
Fredholm theory, and even more surprisingly, it is not much more
difficult than the classical one. This leads to our theory of
M-polyfolds (or more generally polyfolds) and an adapted Fredholm
theory. There are even finite-dimensional polyfolds and Figure
\ref{Fig1.0} shows a finite-dimensional M-polyfold having a
one-dimensional "submanifold". In this new theory we can formalize
new structures which could not be formalized before. This gives a
unified perspective on a variety of theories in symplectic geometry.
It also seems that the theory should have applications in other
fields as well, since the addressed analytical issues arise in
geometric pde's of Riemannian geometry as well as the theory of
nonlinear pde in general. In this note we will explain some parts of
the theory. For the proofs and further discussions we refer the
reader to \cite{HWZ1,HWZ2,HWZ3,EH} and \cite{EGH}.

Acknowledgement: These notes are an expanded version of a talk given
at the CDM 2004 Conference at Harvard, November 2004. I would like
to thank C. Abbas, F. Bourgeois, K. Cieliebak, O. Cornea, T. Ekholm,
Y. Eliashberg,  A. Givental, F. Lalonde, J. Latschev, D. McDuff, K.
Mohnke,  D. Salamon, F. Schlenk, M. Schwarz, K. Wehrheim, K.
Wysocki, and E. Zehnder for many stimulating discussions. The
writing of these notes profited from a workshop in Montreal in 2004
and the SFT-workshop in Leipzig, May 2005, in which the present
theory was discussed in great detail. My thanks go to the DFG and
NSF for their support. The current material was presented with a lot
of detail in a two semester course at the Courant Institute and I
profited from many discussions with B. Bramham, M. Chance, J.
Coffey, E. Dupont, J. Fish, B. Gurel, Z. Han, R. Hind, U.
Hryniewicz, S. Lisi, B. Madero, A.S. Momin, S. Pinnamaneni, A.
Savelyev, R. Siefring, and C. Wendl.

The organization of this paper is as follows:

 \tableofcontents

\section{New Smoothness Concepts and  Spaces}
In a first step we introduce a new concept of smoothness for a
Banach space and define the notion of a smooth map in this new
context. One might, alternatively, view our first definition as a
new interpretation for certain classes of interpolation
spaces\footnote{Interpolation theory is an important part of
functional analysis and the theory of function spaces. We refer to
the comprehensive book by Triebel \cite{Tr}}.
\subsection{Concepts of Smoothness}
\begin{definition}
Let $E$ be a Banach space. A sc-smooth structure on $E$ is given by
a nested sequence of Banach spaces $E_m$, $m\in {\mathbb N}$,
satisfying
\begin{itemize}
\item[1)]  For $m\leq n$ the space $E_n$ is a linear subspace of
$E_m$ and $E_0=E$.

\item[2)] The inclusion $E_n\rightarrow E_m$ for $m<n$ is a
compact operator.

\item[3)] The vector space $E_{\infty}$ defined by
$$
E_{\infty}=\bigcap_{m\in{\mathbb N}} E_m
$$
is dense in every $E_m$.
\end{itemize}
\end{definition}
{ Of course, on can build a linear functional analytic theory on
Banach spaces with sc-structures. There is a large body of such a
theory usually as part of interpolation theory, where the focus is
quite different from ours. The interpretation of a scale as a
generalization of a smooth structure as being developed below
seems to be new.}

Here is an important example.
\begin{example}\label{ex1}
Let $0<\delta_0<\delta_1<...$ be a strictly increasing sequence of
weights. We denote by $E$ the Banach space consisting of maps
$u:{\mathbb R}\times S^1\rightarrow {\mathbb R}^{N}$ of Sobolev
class $H^3_{loc}$ so that for every multi-index $\alpha$ of order at
most $3$ the weighted partial derivative
$$
(s,t)\rightarrow |D^\alpha u(s,t)|\cdot e^{\delta_0|s|}
$$
belongs to $L^2$. We define the sc-structure on $E$ by taking $E_m$
to consist of maps of regularity $(m+3,\delta_m)$, i.e. $m+3$
derivatives integrable with the weight associated to $\delta_m$.
That $E_m$ defines a sc-structure on $E$ follows from the compact
Sobolev embedding theorem for bounded domains and the fact that the
weights are strictly increasing.
\end{example}

If $U\subset E$ is an open subset we define a sc-smooth structure to
be the nested sequence $U_m=U\cap E_m$. Given a sc-smooth structure
on $U$ we observe that $U_m$ inherits a sc-smooth structure by
defining ${(U_m)}_k=U_{m+k}$. We shall write $U^m$ for the sc-space
defined by
$$
(U^m)_k= U_{m+k}.
$$
Given two sc-spaces $E$ and $F$ there is a well-defined direct sum
$E\oplus F$. We would like to note here that there are three
important linear concepts. The first is that of a linear sc-operator
$T:E\rightarrow F$, which by definition is a linear operator
inducing bounded operators between the same levels. The next one is
that of a linear sc-Fredholm operator (to be explained later), and
the latter is that of a sc$^+$-operator $A:E\rightarrow F$. By
definition this is a sc-operator inducing one from $E$ to $F^1$,
i.e. $E_m\rightarrow F_{m+1}$ for every level $m$.

If $U\subset E$ is open we define its tangent as $TU=U^1\oplus E$.
In particular
$$
(TU)_m = U_{m+1}\oplus E_m.
$$
A map $f:U\rightarrow V$, where $U$ and $V$ are open in sc-Banach
spaces, is said to be sc$^0$ provided it induces a continuous map
between every level. Next we define the notion of a sc-smooth map.
We give two equivalent definitions.
\begin{definition}
Let $E$ and $F$ be sc-Banach spaces and $U\subset E$ an open subset.
The sc$^0$-map $f:U\rightarrow F$ is said to be sc$^1$ if the
following holds:
\begin{itemize}
\item For every $x\in U_1$ there exists a bounded linear operator
$Df(x)\in L(E_0,F_0)$ so that for $h\in E_1$
$$
\lim_{h\rightarrow 0}\frac{1}{\parallel h\parallel_1} \parallel
f(x+h)-f(x)-Df(x)h\parallel_0=0.
$$
\item The map $Tf:TU\rightarrow TF$ defined by
$$
(Tf)(x,h)=(f(x),Df(x)h)
$$
is of class sc$^0$.
\end{itemize}
\end{definition}
Let us observe that for every $x\in U_1$ there can be at most one
map $Df(x)$ with the properties described above. This map will be
called the linearization of $f$ at $x$. We call $Tf:TU\rightarrow
TF$ the tangent map of the sc$^1$-map $f:U\rightarrow F$. Let us
observe that $Tf$ for every $m\in {\mathbb N}$ is given by
$$
Tf:U_{1+m}\oplus E_m\rightarrow F_{m+1}\oplus F_m:(x,h)\rightarrow
(f(x),Df(x)h).
$$
There is an equivalent definition  for being sc$^1$ which relates it
to the notion of $C^1$-map between different levels.

\begin{definition}
Let $E$ and $F$ be sc-smooth Banach spaces and $U\subset E$ an
open subset. A sc$^0$-map $f:U\rightarrow F$ is said to be sc$^1$
provided the following holds:
\begin{itemize}\label{sc-1}
\item[1)] For every $m\geq 1$ the induced map
$$
f:U_m\rightarrow F_{m-1}
$$
is of class $C^1$. In particular the derivative gives the
continuous map
$$
U_m\rightarrow L(E_m,F_{m-1}):x\rightarrow Df(x).
$$
\item[2)] For $x\in U_m$ and $m\geq 1$ the map $Df(x)$ induces
 a continuous linear operator
$Df(x):E_{m-1}\rightarrow F_{m-1}$ and the resulting map
$$
U_m\times E_{m-1}\rightarrow F_{m-1}:(x,h)\rightarrow Df(x)h
$$
is continuous.
\end{itemize}
\end{definition}
As already emphasized
\begin{proposition}
The two definitions for a sc$^0$-map $f:U\rightarrow F$ to be of
class sc$^1$ are equivalent.
\end{proposition}
For a proof see \cite{HWZ1}. A sc$^1$-map $f:U\rightarrow V$ has a
well-defined tangent map
$$
Tf:TU\rightarrow TV
$$
and inductively we can define the notion of being sc$^k$. An
important result is the validity of the chain rule:
\begin{theorem}[Chain Rule]\ \
If $f:U\rightarrow V$ and $g:V\rightarrow W$ are sc$^1$ so is
$g\circ f$ and $T(g\circ f)=(Tg)\circ (Tf)$.
\end{theorem}
In view of the second characterization of being sc$^1$ it is not
clear at all that a chain rule has to hold. In fact, as the proof
reveals, it just works.

 We give an example of a sc$^1$-map which
will be important in the application of the theory to SFT.
\begin{example} Recall the
sc-space $E$ of maps ${\mathbb R}\times S^1\rightarrow {\mathbb
R}^N$ from Example \ref{ex1}. We have an action of the group
${\mathbb R}\times S^1$ by sc-operators on $E$ defined by
$$
((c,\rho)\ast u)(s,t)=u(s+c,t+\rho),
$$
where $c\in {\mathbb R}$ and $\rho\in {\mathbb R}/{\mathbb Z}=S^1$.
The map
$$
{\mathbb R}\times S^1\times E_m\rightarrow
E_m:((c,\rho),u)\rightarrow (c,\rho)\ast u
$$
is continuous for every $m$, so that it defines a sc$^0$-map
$$
\Phi:({\mathbb R}\times S^1)\times E\rightarrow E.
$$
The important fact is now that the map $\Phi$ is sc$^\infty$. The
proof is somewhat lengthy and a variation of the proof can be found
in \cite{HWZ1}.
\end{example}

{ At this point we can develop a whole theory of manifolds build on
the pseudogroup of sc-diffeomorphisms. The fact that the whole
manifold theory and its constructions are functorial allow to build
a parallel theory based on this pseudogroup. With other words we
could define a second countable Hausdorff space $X$ to have a
sc-smooth structure provided it is equipped with an atlas so that
the transition maps are sc-smooth. Note that $X$ will inherit a
filtration $X_m$. Also $X$ has a tangent bundle $TX\rightarrow
X^1$.}

Here is an important example.
\begin{example}
Let $M$ be a complete Riemannian manifold and $\Phi:M\rightarrow
{\mathbb R}$ a Morse-function. Let us assume for simplicity that the
critical points can be totally ordered by $a < b$ via $f(a)<f(b)$.
For every $a$ fix a sequence $\delta^a=$ of weights
$\delta_0^a=0<\delta_1^a<..$. Usually the limit should be finite and
smaller than the spectral gap of the Hessian around $a$. Then denote
for $a<b$ by $X(a,b)$ the quotient of $H^2$-maps (Sobolev class)
connecting at $-\infty$ the point $a$ with $b$ at $+\infty$ by the
obvious ${\mathbb R}$-action. The space has a sc-smooth structure
where level-$m$-elements are represented by
$(H^{2+m,\delta_m^a,\delta_m^b})$-maps, i.e. we have different
exponential weights at $\pm\infty$. The local model is the
codimension one sc-subspace $E$ of $H^2({\mathbb R},{\mathbb R}^n)$
consisting of functions $h=(h_1,..,h_n)$ with $h_1(0)=0$, equipped
with the sc-structure $E_m$ given by
$E_m=H^{m+2,\delta_m^a,\delta_m^b}\cap E$. We explain this with the
special case $M={\mathbb R}^N$ and also assume that the weight
sequence is independent of the critical point. Let $a\neq b$ be two
different points in ${\mathbb R}^N$. Pick a smooth map
$\varphi:{\mathbb R}\rightarrow {\mathbb R}^N$ so that
$\varphi(s)=a$ for $s<<0$ and $\varphi(s)=b$ for $s>>0$. Define
$\hat{X}=\varphi+H^2({\mathbb R},{\mathbb R}^N)$. Then we have by
time-shift a natural ${\mathbb R}$-action on $\hat{X}$. Denote the
quotient $\hat{X}/{\mathbb R}$ by $X(a,b)$ with its induced quotient
topology. The space $H^2({\mathbb R},{\mathbb R}^N)$ has a
sc-structure where $E_m$ consists of maps of regularity
$(m+2,\delta_m)$. Given a smooth representative $u$ for a class
$[u]\in X(a,b)$, i.e. $u=\varphi+v$ with $v\in E_\infty$, we can
define the inverse of a chart as follows. Since $u$ connects two
different points $a$ and $b$ there exists a times $t_0$ with
$u'(t_0)\neq 0$. We may assume without loss of generality that
$t_0=0$ by replacing $u$ by $t_0\ast u$. Let $\Sigma$ be the
hyperplane orthogonal to $u'(0)$. Then define a codimension one
subspace of $E$ to consist of all $h$ with $h(0)\in\Sigma$. One can
show that $H$ has a one-dimensional sc-complement in $E$. For $h\in
H$ in a $H^2$-neighborhood of $0$ the map
$$
h\rightarrow [u+h]
$$
is a homeomorphism onto an open neighborhood of $[u]\in X(a,b)$. One
can show that the collection of the inverses of all these maps
defines an atlas of charts with sc-smooth transition maps. Hence
$X(a,b)$ carries the structure of a sc-manifold.  Going back to our
original case of maps into $M$ we can modify the construction for
the ${\mathbb R}^N$-case using the exponential map for a suitable
Riemannian metric on $M$ and can construct a sc-manifold structure
in the general case as well.

\end{example}

\subsection{M-Polyfolds}
The local models for the new spaces, which are needed to construct
ambient spaces for the moduli spaces (for example occurring in SFT),
can now easily be constructed ( at least those which have a
"manifold flavor" ).  We call these new type of spaces polyfolds, We
begin with M-polyfolds, where the $M$ stands for "manifold flavor".
Polyfolds, which are defined later have an orbifold flavor. These
spaces have locally varying dimensions. We also exhibit a (rather
tame) finite-dimensional example.
\subsubsection{M-Polyfolds}
We will describe first M-polyfolds which carry a manifold flavor
in contrast to { polyfolds which are in general potentially
complicated objects resembling something like "orbifolds with
varying dimensions and with boundary with corners"}. Let us call a
subset $C$ of some finite-dimensional vector space $A$ a partial
cone if there is a linear isomorphism $T:A\rightarrow {\mathbb
R}^n$ mapping $C$ onto $[0,\infty)^k\times {\mathbb R}^{n-k}$.

\begin{definition}
Let $V$ be a (relatively) open subset of some partial cone $C$,
$E$ a Banach space with a sc-smooth structure and
$\pi_v:E\rightarrow E$ a family of sc-projections so that the
induced map
$$
V\oplus E\rightarrow E:(v,e)\rightarrow \pi_v(e)
$$
is sc-smooth. Then we call the triple ${\mathcal S}=(\pi,E,V)$ a
sc-smooth splicing.
\end{definition}
Every splicing ${\mathcal S}=(\pi,E,V)$ is accompanied by a
complementary splicing ${\mathcal S}^c=(Id-\pi,E,V)$. Observe that a
splicing ${\mathcal S}$ decomposes the space $V\oplus E$ as a
fibered sum over $V$, namely a point $(v,x)$ can we decomposed as
$$
(v,x) = (v,u_v +u^c_v),
$$
where $\pi_v(u)=u$ and $\pi_v(u^c)=0$. Also note that the
sc-smoothness of $(v,e)\rightarrow \pi_v(e)$ is a rather weak
condition. In fact the dimension of the image of $\pi_v$ is in
general locally not constant.

The map $\Phi:V\oplus E\rightarrow E:(v,e)\rightarrow \pi_v(e)$ is
sc-smooth. Taking its tangent map we can define
$$
\Pi_{(v,\delta v)}:TE\rightarrow TE:(e,\delta e)\rightarrow
(\Phi(v,e),D\Phi(v,e)(\delta v,\delta e))
$$
which has the property that the induced map
$$
TV\oplus TE\rightarrow TE:(a,b)\rightarrow \Pi_a(b)
$$
is sc-smooth since it is modulo the identification $TV\oplus
TE=T(V\oplus E)$ the tangent map of $\Phi$. One easily verifies
that $(\Pi,TV,TE)$ defines a sc-smooth splicing. We call it the
tangent of the splicing ${\mathcal S}$ and denote it by
$T{\mathcal S}$.

Let ${\mathcal S}=(\pi,E,V)$ be a sc-smooth splicing. Then the
associated splicing core is the subset $K=K^{\mathcal S}$ of
$V\oplus E$ consisting of all pairs $(v,e)$ with $\pi_v(e)=e$.
Observe that we have a natural map
$$
K^{T{\mathcal S}}\rightarrow K^{\mathcal S}:(v,\delta v,e,\delta
e)\rightarrow (v,e).
$$
Clearly the fiber over any point is a sc-Banach space in a natural
way. We define the tangent of $K^{\mathcal S}$, denoted by
$TK^{\mathcal S}$, by
$$
TK^{\mathcal S}:=K^{T{\mathcal S}}.
$$
In fact it is useful to keep track of the underlying splicings and
the above should be read
$$
T(K^{\mathcal S},{\mathcal S})= (K^{T{\mathcal S}},T{\mathcal S}).
$$
Now we can define our new local models for spaces.
\begin{definition}
A local M-polyfold model consists of a pair $(O,{\mathcal S})$ where
$O$ is an open subset of the splicing core $K^{\mathcal S}$
associated to the sc-smooth splicing ${\mathcal S}$. The tangent
$T(O,{\mathcal S})$ of the local M-polyfold model $(O,{\mathcal S})$
is defined by
$$
T(O,{\mathcal S}) = (K^{T{\mathcal S}}|O,T{\mathcal S}),
$$
where $K^{T{\mathcal S}}|O$ denotes the collection of all points
in $K^{T{\mathcal S}}$ which project under the canonical
projection
$$
K^{T{\mathcal S}}\rightarrow K^{\mathcal S}
$$
onto $O$.
\end{definition}
The above discussion gives us  a natural projection
$$
T(O,{\mathcal S})\rightarrow (O,{\mathcal S}):(v,\delta v,e,\delta
e)\rightarrow (v,e).
$$
In the following we shall write $O$ instead of $(O,{\mathcal S})$,
but observe that ${\mathcal S}$ is part of the structure. Note
that for an open subset $O$ of a splicing core we have an induced
filtration. Hence we may talk about sc$^0$-maps.  We continue by
introducing the notion of a sc$^1$-map between open sets of
splicing cores.

\begin{definition}
Let $O$ and $O'$ be open subsets of splicing cores. Assume that
$f:O\rightarrow O'$ is a sc$^0$-map. We say that $f$ is sc$^1$
provided the map
$$
(v,e)\rightarrow f(v,\pi_v(e))
$$
which is defined on some open subset of $C\oplus E$  and takes
image in $O'$ (which we view lying in the obvious  Banach space
with sc-smooth structure) is sc$^1$.
\end{definition}

A sc$^1$-map induces in a canonical way  a tangent map
$$
Tf:TO\rightarrow TO'.
$$
To see this start with the map
$$
\hat{f}:(v,e)\rightarrow f(v,\pi_v(e)) =
(f_1(v,\pi_v(e)),f_2(v,\pi_v(e)))=(\hat{f}_1(v,e),\hat{f}_2(v,e))
$$
which is defined on an open subset $\hat{O}$ of a sc-smooth Banach
space and which takes its image in a sc-smooth Banach space $G$. By
assumption it is sc$^1$. On easily verifies that the map $Tf$
defined by
$$
Tf(v,\delta v,e,\delta e) = (T\hat{f}_1(v,e,\delta v,\delta
e),T\hat{f}_2(v,e,\delta v,\delta e))
$$
maps $K^{{T\mathcal S}}|O$ into $K^{T{\mathcal S}'}|O'$. This is
by definition the induced tangent map. We have
\begin{theorem}[Chain Rule for sc$^1$-maps]\ \
Let $O,O',O''$ be open sets in splicing cores and $f:O\rightarrow
O'$ and $g:O'\rightarrow O''$ be sc$^1$. Then $g\circ f$ is sc$^1$
and $T(g\circ f) =Tg\circ Tf$. Moreover $Tf$ and $Tg$ are sc$^0$.
\end{theorem}
This is a consequence of the sc-chain rule, the definition and the
fact that our reordering of the terms in our definition of the
tangent map is consistent. Hence given a sc$^1$-map $f:O\rightarrow
O'$ between open sets of splicing cores we obtain an induced tangent
map $Tf:TO\rightarrow TO'$. Inductively we can define the notion of
being sc$^k$.

Here one should point out the following. If we recall that the
constructions of Differential Geometry are functorial with respect
to the input being: a) the notion of a smooth map between two open
subsets of Euclidean spaces and b) the chain rule and the
functoriality of the tangent functor, then we can easily imagine if
we replace open sets in Euclidean spaces by open sets in splicing
cores and smooth maps by sc-smooth maps that many constructions of
Differential Geometry carry over and that many more constructions
become possible. We should however note that in finite dimensions
the existence of a smooth partition of unity is automatic, whereas
in general (with the exception of sc-Hilbert spaces) it has to be
required\footnote{For example it is know that $L^p$-spaces with
$1<p<\infty$ and $p\neq 2$ do not allow smooth partitions of unity.
The best possible differentiability $k(p)$ is a monotonic function
in $p$ with $k(p)\rightarrow \infty$ for $p\rightarrow\infty$.}. The
situation is, however, somewhat better in the sc-situation than in
the Banach space case in view of the following criterion for
sc-smoothness of a real-valued function, which is proved in
\cite{HWZ1}.

\begin{proposition}
Let $E$ be a Banach space with a sc-smooth structure and let
$U\subset E$ be an open subset. Assume that $f:U\rightarrow{\mathbb
R}$ is sc-continuous and that the induced maps
$$f_m:=f\vert_{U_m}:U_m\to {\mathbb R},\quad m\geq 0,$$
are of class $C^{m+1}$. Then $f$ is of class sc$^{\infty}$.
\end{proposition}

This should allow to define sc-smooth partitions on spaces which
classically do not admit smooth partitions of unity. It would be
interesting to see some worked out examples for interesting (i.e.
relevant for applications) sc-structures.

Armed with the philosophical point of view that most constructions
of Differential Geometry should carry over if we replace open sets
of Banach spaces by open sets of splicing cores, we can introduce
the notion of a M-polyfold. This is in the new context the object
corresponding to the classical notion of a manifold.
\begin{definition}
Let $X$ be a second countable Hausdorff space. A M-polyfold chart
is a triple $(U,\varphi,{\mathcal S})$, where $U$ is an open
subset of $X$ and $\varphi:U\rightarrow K^{\mathcal S}$ a
homeomorphism onto an open subset of a splicing core. We say two
charts are compatible if the transition map between open subsets
of splicing cores is sc-smooth in the sense defined above. A
maximal atlas of sc-smoothly compatible M-polyfold charts is
called a M-polyfold structure on $X$.
\end{definition}
Let us observe that a M-polyfold is necessarily metrizable. If $X$
is a M-polyfold so is $X^1$ and moreover $X^n$ for any $n\geq 1$.
Here $X^n={(X^{n-1})}^1$.  Given a M-polyfold $X$ we can construct
its tangent $TX$ in a natural way. The projection
$$
TX\rightarrow X^1
$$
is a sc-smooth map. One may view $TX$ as a bundle over $X^1$. As it
will turn out we need to introduce the notion of a strong bundle in
order to develop a Fredholm theory. The tangent bundle will in
general not be a strong bundle.

For the convenience of the reader let us give some examples
illustrating the new notions. We begin with an example for a
splicing.
\begin{example}
Let $E=L^2({\mathbb R})$ be equipped with the sc-structure defined
by $E_m=H^{m,\delta_m}$, i.e. maps of Sobolev class $H^m_{loc}$ with
derivatives up to order $m$ weighted by $e^{\delta_m |s|}$ belonging
to $L^2$. Here $\delta_0=0<\delta_1<..$ is a strictly increasing
sequence of weights. Pick a smooth compactly supported map
$\gamma:{\mathbb R}\rightarrow [0,1]$ with
$$
\int_{\mathbb R} \gamma(t)^2 dt =1.
$$
Next we put $V={\mathbb R}$ and define for $t\in V$ a family of
sc-projections $\pi_t$ by $\pi_t=0$ for $t\leq 0$ and for $t>0$
$$
\pi_t(u) =\langle u,\gamma_t\rangle\cdot \gamma_t,
$$
where
$$
\gamma_t(s) = \gamma(s+e^{\frac{1}{t}}).
$$
One can show that ${\mathcal S}:=(\pi,E,V)$ is a sc-smooth splicing.
The splicing core $K^{\mathcal S}$ is homoeomorphic to
$$
X = \left((-\infty,0]\times \{0\}\right) \bigcup
\left((0,\infty)\times {\mathbb R}\right).
$$
From this it follows that $X$ can be equipped with a M-polyfold
structure. The reader will easily modify this example to construct
splicings where the local dimensions vary between $1$ and any given
natural number $N$. As a consequence one can show that the subspace
$X$ of ${\mathbb R}^3$ defined below admits a M-polyfold structure:
\begin{eqnarray*}
X&=&\{(x,0,0)\ |\ x\in(-\infty,-1]\cup [1,2]\}\bigcup \{(x,y,0)\ |\
x^2+y^2<1\}\\
&&\bigcup \{(x,y,z)\ |\ |x-3|^2+y^2+z^2 <1\}.
\end{eqnarray*}
Then $X^m$ is independent of $m$. Moreover $TX$ would be over the
open unit disk a plane bundle and over the line a line bundle, etc.
There are, in fact, much more complicated examples with the
dimensions even allowed to locally vary between finite and infinite.
\end{example}

An interesting feature is that sc-smooth recognize corners. To prove
this requires some efforts and we refer the reader for a proof in
\cite{HWZ1}. As a consequence a M-polyfold has its corner structure
as an invariant. Take a M-polyfold chart $U\rightarrow O$, then $O$
is in an open subset of a splicing core $K=\{(v,e)\in V\oplus E\ |\
\pi_v(e)=e\}$. Here $V$ is open in a partial cone, say
$[0,\infty)^k\times {\mathbb R}^{n-k}$. We can assign to a point
$x\in X$ the number $d(x)$ of vanishing first $k$-coordinates. It
turns out that this does not depend on the choice of local
coordinates and every point $x$ has an open neighborhood so that
$d|U(x)\leq d(x)$.
\begin{definition}
For a M-polyfold $X$ the map $d:X\rightarrow {\mathbb N}$ is called
the degeneration map.
\end{definition}
This map  will be important in our Fredholm theory with operations.
We need another definition.
\begin{definition}
Given a M-polyfold $X$ we call the closure of a connected component
$F$ of $X(1)=\{x\in X\ |\ d(x)=1\}$ a face.
\end{definition}
It is an easily established fact that around every point $x_0\in X$
there exists an open neighborhood $U=U(x_0)$ so that every $x\in U$
belongs to exactly $d(x)$ many faces of $U$. Globally it is always
true that $x\in X$ belongs to at most $d(x)$ many faces and strict
inequality is possible. For example a two-dimensional closed domain
with one corner point, homeomorphic to the closed disk is not
face-structured.

\begin{definition}
We call a M-polyfold face-structured if every point $x$ belongs to
$d(x)$ many faces.
\end{definition}

Face-structure M-polyfolds and polyfolds will be important in SFT,
or more generally in a Fredholm theory with operations, since they
have an interesting algebraic structure. For known facts about
finite-dimensional manifolds with boundary and corners see for
example \cite{Janich} and \cite{Laures}.

We continue with our Morse-Theory example in a separate subsection.

\subsubsection{Example of a M-polyfold in Morse-theory}
 Consider again our
Morse function $\Phi:M\rightarrow {\mathbb R}$. For any finite
sequence $(a_0,..,a_k)$ of critical points with $a_i<a_{i+1}$ and
$k\geq 1$ define $X(a_0,..,a_k)=X(a_0,a_1)\times..\times
X(a_{k-1},a_k)$ and let $\overline{X}$ be the disjoint union of all
these $X(a_0,..,a_k)$. One can equip $X$ with a natural second
countable Hausdorff topology inducing on all parts the already
defined topology and in addition the closure of $X(a,b)$ contains
all $X(a_0,...,a_k)$ where $a_0=a$ and $a_k=b$. Moreover there is a
natural M-polyfold structure on $\overline{X}$ so that
$\overline{X}$ is faced structured and $d(x)=k-1$ provided $x\in
X(a_0,..,a_k)$ for some sequence $(a_0,..,a_k)$ of ordered critical
points. It requires some work to write down the relevant splicing
cores and the charts. The main point  is, of course, the
understanding of the space $\overline{X}$ near a broken trajectory.
For this we have to introduce a particular splicing which we discuss
now. Similar versions will be crucial for the constructions of
polyfolds in SFT. We will only discuss a model situation and refer
the reader to \cite{HWZ1} for full details.

Let us assume we are given three mutually different points $a,b$ and
$c$ in ${\mathbb R}^N$. We have seen that $X(a,c)$ is a sc-manifold
in a natural way given a sequence of weights $(\delta_i)$ starting
with $\delta_0=0$. Assuming that $a$ and $c$ are critical points for
a Morse function $\Phi$ the gradient lines connecting $a$ with $c$
(modulo parametrization) would lie in $X(a,c)$. The space of these
gradient lines is in general not compact and a gradient line might
split into a broken gradient line first going from $a$ to another
critical point $b$ and then to $c$. In order to compactify the space
of gradient lines (which we might view as the solution space of a
nonlinear  elliptic problem ) we have to add suitable broken ones.
If we want to develop a Fredholm theory for which the compactified
space is the solution space, the ambient space needs to contain
broken trajectories. A natural choice as a set in our case is
obviously
$$
\bar{X}(a,c)=\left(X(a,b)\times X(b,c)\right) \bigcup X(a,c).
$$
The spaces $X(a,b), X(b,c)$ and $X(a,c)$ have natural paracompact
second countable topologies. One can show that there is a natural
second countable paracompact topology on $\bar{X}(a,c)$ inducing on
$X(a,c)$ and $X(a,b)\times X(b,c)$ the given topology so that
$X(a,c)$ is dense in $\bar{X}(a,c)$. The topological space
$\bar{X}(a,c)$ is not obviously homeomorphic to any open subset of a
Banach space. However, it is  a nontrivial fact that it is
homeomorphic to an open subset of some splicing core. Moreover, all
these local homeomorphisms can be picked in such a way that the
transition maps are sc-smooth.

The construction of the relevant splicing is closely related to some
gluing construction. We begin with "nonlinear gluing" which quite
often in literature is referred to as pre-gluing.  It associates to
a curve connecting $a$ with $b$, and one connecting $b$ with $c$,
and a gluing parameter $r\in (0,1)$ a curve connecting $a$ with $c$.
Then we define a (linear) gluing and anti-gluing for vector fields
along the underlying given curves. Finally we will show how these
constructions are related to each other and how we can construct an
associated splicing. It will be a punch-line that the splicing idea
can be viewed as a generalization of some constructions arising
around the gluing procedure.

Let $\beta:{\mathbb R}\rightarrow [0,1]$ be a smooth cut-off
function so that
\begin{eqnarray*}
&\beta(s)=1\ \hbox{ for} \ s \leq -1&\\
&\beta'(s)<0\ \hbox{ for} \ s\in (-1,1)&\\
&\beta(s)+\beta(-s)=1\ \hbox{ for all} \  s\in {\mathbb R}.&
\end{eqnarray*}
Let us assume we are given $u,v :{\mathbb R}\rightarrow {\mathbb
R}^N$ with $u(-\infty)=a$, $u(+\infty)=b$, $v(-\infty)=b$ and
$v(+\infty)=c$. These two maps are representatives of classes
$[u]\in X(a,b)$ and $[v]\in X(b,c)$. For a real number $R\geq 0$ we
define the glued map $\oplus_R(u,v)$ by
$$
\oplus_R(u,v)(s) = \beta(s-\frac{R}{2})u(s)
+(1-\beta(s-\frac{R}{2}))v(s-R).
$$
Of interest for us will be the class $[\oplus_R(u,v)]\in X(a,c)$. We
will also define a gluing for $R=\infty$ by
$$
\oplus_\infty(u,v)=(u,v).
$$
The number $R\in [0,\infty)\cup \{\infty\}=:[0,\infty]$ we will call
the gluing length. At this point we have defined a gluing for any
gluing length in $[0,\infty]$. On the level of equivalence classes
we would like to view $([u],[v])$ as on the boundary of the family
$\{[\oplus_R(u,v)]\ |\ R\in [0,\infty)\}$. For this we have to
identify the family, say, with $[0,1)$ so that $([u],[v])$
corresponds to $1$. One has to be precise here, since there are many
ways of identifying $[0,1]$ with $[0,\infty]$. Also a precise choice
is required for the definition of the M-polyfold structure
(sc-smoothness of the transition maps!). Let us call a number $r\in
[0,1]$ a gluing parameter. In order to make a consistent
construction which leads to sc-smooth transition maps only the
following piece of data, namely a gluing profile, is needed.
\begin{definition}
A gluing profile $\varphi$ is a diffeomorphism
$$
\varphi:(0,1]\rightarrow [0,\infty)
$$
\end{definition}
A gluing profile is obviously a rule how a gluing parameter is
converted into a gluing length. Clearly $r=0$ is associated to
$R=\infty$. We refer the reader to \cite{EHWZ1}, where different
gluing profiles are studied in the context of Deligne-Mumford Theory
of stable Riemann surfaces. Usually, in any application a gluing
profile has to satisfy certain growth conditions. A very useful
gluing profile is the "exponential profile"
$$
\varphi(r) = e^{\frac{1}{r}}-e
$$
and we will use it in the following. Using the letters $r$ and $R$
one should always have in mind that $R=\varphi(r)$.

\begin{figure}[htbp]
\mbox{}\\[1ex]
\centerline{\relabelbox
\epsfxsize 4.4truein \epsfbox{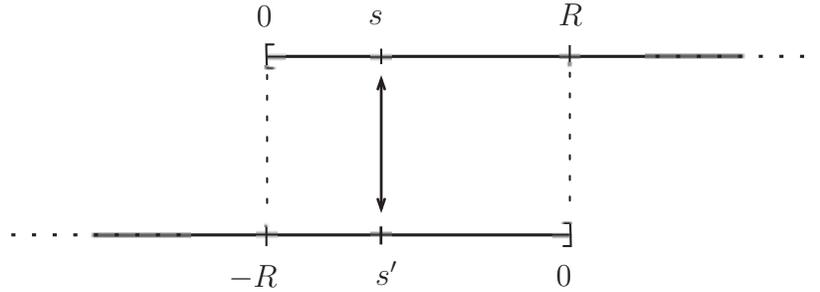}
\relabel {z}{$0$}
\relabel {az}{$0$}
\relabel {s}{$s$}
\relabel {sp}{$s'$}
\relabel {mr}{$-R$}
\relabel {r}{$R$}
\endrelabelbox}
\caption{Identification $-R+s=s'$}\label{Figure3}
\mbox{}\\[1ex]
\end{figure}

Let us recall that the inverses of the charts for $X(a,b)$ were of
the form
$$
h\rightarrow [u+h],
$$
where $h(0)\in \Sigma$ and $\Sigma$ is a hypersurface in ${\mathbb
R}^N$. One should interpret $u+h$ as $\exp_u(h)$, where $h$ is a
vector field along the underlying curve $u$. Similarly
$$
k\rightarrow [v+k].
$$
Using the gluing construction we can define
$$
(r,h,k)\rightarrow [\oplus_R(u+h,v+k)].
$$
Any curve $[w]$ near to the broken one $([u],[v])$ can be written in
such way. Of course there is a large ambiguity. The class $[w]$ can
be written  in many such ways (taking different choices of $u$ and
$v$) and in the following we introduce the concept of anti-gluing to
get rid of this ambiguity.

We begin by introducing the "linear gluing" on the level of vector
fields. We use the same formula as before
$$
\oplus_R(h,k)(s)=\beta(s-\frac{R}{2})h(s)
+(1-\beta(s-\frac{R}{2}))k(s-R).
$$
It is, of course, important to find a proper interpretation for
$\oplus_R(h,k)$, particularly since we want to generalize our model
situation to the Morse-theory situation for a manifold $M$. The
right interpretation for $\oplus_R(h,k)$ is that of a vector field
along $\oplus_R(u,v)$. Then the obvious relationship
$$
\oplus_R(u+h,v+k) = \oplus_R(u,v) +\oplus_R(h,k)
$$
can be rewritten as
$$
\oplus_R(\exp_u(h),\exp_v(k))= \exp_{\oplus_R(u,v)}(\oplus_R(h,k)).
$$
One can define gluing in the general case in such a way that this
formula is true for certain Riemannian metrics. After having defined
the "linear gluing" on the level of vector fields along the
underlying curves we introduce now the "linear anti-gluing"
$\ominus_R(h,k)$. Again $h$ and $k$ are vector fields along $u$ and
$v$ as described before. Then we define
$$
\ominus_R(h,k)(s)=-(1-\beta(s-\frac{R}{2}))h(s)+\beta(s-\frac{R}{2})k(s-R).
$$
Observe that the map
$$
(h,k)\rightarrow (\oplus_R(h,k),\ominus_R(h,k))
$$
is an isomorphism (for fixed $R$). The right interpretation for
$\ominus_R(h,k)$ is that of a map
$$
{\mathbb R}\rightarrow T_b({\mathbb R}^N).
$$
With the gluing constructions at hand we can define the splicing.
Denote by $E$ the sc-space consisting of pairs $(h,k)$ of vector
fields along $u$ and $v$, respectively, which satisfy
$h(0)\in\Sigma$ and $k(0)\in\Sigma'$.

Then $E$ has for every $r\in [0,1)$ two distinguished
sc-complementary subspaces, namely $\ker(\ominus_R)$ and
$\ker(\oplus_R)$:
$$
E=\ker(\ominus_R)\oplus_{sc}\ker(\oplus_R).
$$
In case $r=0$ we have $\ker(\ominus_\infty)=E$ and
$\ker(\oplus_\infty)=\{0\}$. We define $V=[0,1)$ which is an open
subset of the cone $[0,\infty)\subset {\mathbb R}$ and denote by
$\pi_r$ the projection onto $\ker(\ominus_R)$ along
$\ker(\oplus_R)$. It is a nontrivial result that
\begin{theorem}\label{thxz}
The triple $(\pi,E,V)$ is a sc-smooth splicing.
\end{theorem}
Then, which is again a nontrivial result we have
\begin{theorem}
Let $(u,v)$ be a smooth pair of paths connecting $a$ via $b$ to $c$
as described before. Then the map
$$
(r,h,k)\rightarrow [\oplus_R(\exp_u(h),\exp_v(k))]
$$
defined for $(r,h,k)$ sufficiently close to $(0,0,0)$ in the
splicing core $K^{\mathcal S}$ is a homeomorphism onto an open
subset of $\bar{X}(a,c)$. Hence its inverse is a chart with image
being an open subset of the splicing core $K^{\mathcal S}$. Here
${\mathcal S}$ is the splicing core from Theorem \ref{thxz}.
Moreover all these charts together with the charts of the
sc-manifold $X(a,c)$ previously constructed are sc-smoothly
compatible.
\end{theorem}
The previously discussed construction can be brought into a manifold
set-up and we can define a M-polyfold structure on the space of
broken curves connecting a critical point $a$ with $b$. We refer the
reader to \cite{HWZ1} for the precise construction.

\subsubsection{Local Strong sc-Bundles}

Let $E$ and $F$ be two Banach spaces with sc-smooth structures. We
define their $\triangleleft$-product $E\triangleleft F$ which
consists of $E\oplus F$ with the double filtration
$$
(E\triangleleft F)_{m,k} = E_m\oplus F_k
$$
defined for $0\leq k\leq m+1$. Observe that the product is not
symmetric. Given an open subset $U$ of $E$ we define in the obvious
way $U\triangleleft F$. We have a canonical map
$$
U\triangleleft F\rightarrow U.
$$
We refer to the above as a local strong sc-bundle. We define the
$\triangleleft$-tangent space by
$$
T_{\triangleleft}(U\triangleleft F)=(TU)\triangleleft (TF).
$$
Again we have a double filtration with $(m,k)$ and $k\leq m+1$ by
$$
T_{\triangleleft}(U\triangleleft F)_{m,k}= U_{m+1}\oplus E_m\oplus
F_{k+1}\oplus F_k.
$$
Given $U\triangleleft F\rightarrow U$ we can build  the associated
derived sc-spaces
$$
U\oplus F \ \ \hbox{and}\ \  U\oplus F^1.
$$
\begin{definition}
Let $U\triangleleft F\rightarrow U$ and $V\triangleleft G\rightarrow
V$ be two local strong sc-bundles. A sc$_{\triangleleft}^0$-map is a
map
$$
f:U\triangleleft F\rightarrow V\triangleleft G
$$
of the form
$$
f(u,h) = (a(u),\ell(u,h))
$$
inducing sc$^0$-maps between the associated derived sc-spaces, i.e.
we have induced maps
$$
f:U\oplus F\rightarrow V\oplus G
$$
and
$$
f:U\oplus F^1\rightarrow V\oplus G^1
$$
which are sc$^0$.
\end{definition}
We can define a sc$_\triangleleft^1$-notion as follows:
\begin{definition}
We say that the sc$_{\triangleleft}^0$-map $f:U\triangleleft
F\rightarrow V\triangleleft G$ is sc$_\triangleleft^1$ provided it
induces sc$^1$ maps between the associated derived sc-spaces.
\end{definition}
If $f$ is sc$_\triangleleft^1$ we obtain an induced
sc$^0_\triangleleft$-map
$$
T_{\triangleleft}f:T_{\triangleleft}(U\triangleleft F)\rightarrow
T_{\triangleleft}(V\triangleleft G)
$$
defined by
$$
(T_{\triangleleft}f)(u,h,v,b)=(a(u),Da(u)h,\ell(u,v),D\ell(u,v)(h,b)).
$$
The following is easily obtained:
\begin{theorem}[Chain Rule for sc$^1_\triangleleft$-maps]\ \
Let $f:U\triangleleft F\rightarrow R\triangleleft G$ and
$g:V\triangleleft G\rightarrow W\triangleleft H$ be
sc$_\triangleleft^1$, where $U\subset E$ and $V\subset R$ are open
with the induced sc-structures, and $g(U\triangleleft F)\subset
V\triangleleft G$. Then $g\circ f$ and sc$_\triangleleft^1$ and
$$
T_{\triangleleft}(g\circ f)=(T_\triangleleft
g)\circ(T_\triangleleft f).
$$
Moreover $T_{\triangleleft}(g\circ f)$ is sc$_\triangleleft^0$.
\end{theorem}
Inductively we can define the notion of being
sc$^k_\triangleleft$. If $\Phi:U\triangleleft E\rightarrow
V\triangleleft F$ is sc$^k_\triangleleft$ and has the form
$$
\Phi(x,h)=(a(x),\phi(x,h))
$$
and is linear in $h$ we call it a sc$^k_\triangleleft$-vector
bundle map. We call $\Phi$ a sc$_{\triangleleft}$ vector bundle
isomorphism if it is sc$_\triangleleft$-smooth and the same holds
for the inverse.

Given a local strong sc-bundle there are two important classes of
sc-smooth sections. Let $U\triangleleft F\rightarrow U$ be the
bundle. A sc-smooth section $f$ is a map of the form
$$
u\rightarrow (u,\bar{f}(u))
$$
so that the principal part $\bar{f}:U\rightarrow F$ is sc-smooth. A
sc-section $f$ is called a sc$^+$-smooth section provided the
principal part induces a sc-smooth map $\bar{f}:U\rightarrow F^1$.
The tangent $T_\triangleleft f$ of a sc$^+$-section $f$ is a
sc$^+$-section of $TU\triangleleft TF\rightarrow TU$. Let us denote
the space of sc-smooth sections by $\Gamma(U\triangleleft F)$ and
that of sc$^+$-sections by $\Gamma^+(U\triangleleft F)$. Assume that
$\Phi:U\triangleleft F\rightarrow V\triangleleft G$ is a
sc$_{\triangleleft}$-smooth vector bundle isomorphism. Then the
pull-back maps induce isomorphisms
$$
\Gamma(V\triangleleft G)\rightarrow \Gamma(U\triangleleft G)
$$
and
$$
\Gamma^+(V\triangleleft G)\rightarrow \Gamma^+(U\triangleleft G)
$$
The same, of course holds for the push-forward. { As a consequence
of the previous discussion we can construct a theory of strong
sc-bundles over sc-manifolds. More precisely, let $b:Y\rightarrow X$
be a continuous surjective map between second countable Hausdorff
spaces so that the preimage of a point $x\in X$ has the structure of
Banach space. Then we can equip $Y$ with charts preserving the
algebraic structure in the fiber so that the transition maps are
sc$_\triangleleft$-smooth vector bundle isomorphisms. The space $X$
then has the underlying structure of a smooth sc-manifold. Note that
$X$ inherits a filtration $X_m$, whereas $Y$ has a double-filtration
$Y_{m,k}$ with $0\leq k\leq m+1$. We denote by $\Gamma(b)$ the
vector space of sc-smooth sections of $b$ and by $\Gamma^+(b)$ the
vector space of $(+)$-sections.}

Let us explain the philosophy behind this concept again with our
usual Morse-theory situation. We will stay in the sc-manifolds
setting. Let $a<b$ be critical points and consider for $u$
connecting $a$ with $b$ and belonging to $H^2$, sections of class
$H^1$ along $u^{\ast}TM$. Clearly the map
$$
u\rightarrow \dot{u}-\Phi'(u)
$$
maps $u$ of class $H^2$ to a $H^1$-section along $u$. If we insist
on notational grounds that sections preserve the filtration index,
then $H^2$-maps into $M$, lie on the same level as the
$H^1$-sections along them. The following observation is crucial.
Note that it makes sense to talk about $H^2$-sections along an
underlying $H^2$-curve as well, i.e. it makes sense to talk about
sections of a somewhat higher regularity along a base curve then a
priori seems to be needed for the Fredholm theory. Denote by
$Y(a,b)$ the space of equivalence classes defined by the ${\mathbb
R}$-action. We have a canonical map
$$
Y(a,b)\rightarrow X(a,b):[(u,h)]\rightarrow [u].
$$
The fiber over a point $[u]$ has a natural Banach space structure.
One can equip $Y(a,b)$ with a second countable Hausdorff topology so
that the projection map is continuous. Moreover, it possesses the
structure of a strong sc$_\triangleleft$-bundle, where the subspace
$Y_{m,k}$ consists of all $(k+1,\delta^a_k,\delta^b_k)$-sections
along a $(m+2,\delta^a_m,\delta^b_m)$-map. The map
$$
f:[u]\rightarrow [\dot{u}-\Phi'(u)]
$$
defines a sc-smooth section. This section will turn out to be a
Fredholm operator in our generalized sense. Clearly $f$ maps $X_m$
to $Y_{m,m}$. A sc$^+$-section $s$ maps $X_m$ to $Y_{m,m+1}$. It
will turn out that as a consequence of the compactness property for
sc-structures a perturbation of $f$ by $s$ still will be a Fredholm
section. However not for the bundle $b$, but for the bundle with
shifted index $b^1$:
$$
b^1:{(Y^1)}_{m,k}:=Y_{m+1,k+1}\rightarrow {(X^1)}_m:=X_{m+1}.
$$

\subsubsection{M-Polyfold Bundles}

The next step consists in introducing M-polyfold bundles. { For this
we need a particular notion of splicing. Of course, not only the
base should be spliced but also the fiber.}
\begin{definition}
A spliced  sc-fibered Banach scale is a triple ${\mathcal
S}_{\triangleleft}= (\Pi,{E}\triangleleft{H},V)$ where
$\Pi=(\pi,\sigma)$ and $\pi_v:{E}\rightarrow { E}$ and
${\sigma_v}:{H}\rightarrow { H}$ are splicing families parameterized
by $V$.
\end{definition}

Note that the above data gives  splicings ${\mathcal S}_0 =
(\pi,E,V)$ and ${\mathcal S}_1=(\sigma,H,V)$. Then we build the
fibered $\triangleleft$-product
$$
 K^{{\mathcal
S}_\triangleleft}:=K^{{\mathcal S}_0}\triangleleft_V K^{{\mathcal
S}_1}
$$
with the double filtration by $[m,k]$ (We write $[m,k]$ to
indicate that $0\leq k\leq m+1$.) given by
$$
{\left(K^{{\mathcal S}_0}\triangleleft_V K^{{\mathcal
S}_1}\right)}_{m,k} = \left\{(v,e,h)\in V\oplus E_m\oplus H_k\ |\
\pi_v(e)=e,\ \rho_v(h)=h\right\}
$$
The natural projection
$$
(V\oplus E)\triangleleft H\rightarrow V\oplus E
$$
induces a natural projection
$$
K^{{\mathcal S}_{\triangleleft}} \rightarrow K^{{\mathcal S}_0}.
$$
As before we can define a tangent $T{\mathcal S}_\triangleleft$ of
the splicing ${\mathcal S}_\triangleleft$ and associated
$K^{T{\mathcal S}_{\triangleleft}}$ so that $T\pr_1$ induces
(using the definition)
$$
T\pr_1:TK^{{\mathcal S}_{\triangleleft}} \rightarrow TK^{{\mathcal
S}_0}.
$$
We are interested in pairs $(K^{{\mathcal
S}_{\triangleleft}}|O,{\mathcal S}_{\triangleleft})$, where
${\mathcal S}_{\triangleleft}$ is a spliced sc-fibered Banach
scale $(\pi,E\triangleleft H,V)$ and
$$
K^{{\mathcal S}_{\triangleleft}}|O
$$
stands for the preimage under the canonical projection
$$
K^{{\mathcal S}_{\triangleleft}}\rightarrow K^{\mathcal S}
$$
of the open subset of $O$ of $K^{{\mathcal S}_0}$. A
sc$_\triangleleft$-smooth morphism
$$
\Phi:(K^{{\mathcal S}_{\triangleleft}}|{O},{\mathcal
S}_{\triangleleft} )\rightarrow ( {K}^{{\mathcal
S}'_{\triangleleft}}| O',{\mathcal S}'_{\triangleleft})
$$
is a map
$$
{K}^{{\mathcal S}_{\triangleleft}}|{ O} \rightarrow {K}^{{\mathcal
S}'_{\triangleleft}}|{ O}'
$$
of the form
$$
(a,b)\rightarrow (\phi(a),\Phi(a,b)).
$$
so that its composition with the projections is a
sc$_{\triangleleft}$-smooth map. To be more precise if $a=(r,e,h)$
with $r\in V$, $h\in {H}$ and $e\in {E}$, then
$$
(r,e,b)\rightarrow (\phi(r,\pi_r(e)),
\Phi(r,\pi_r(e),\sigma_{r}(b)))
$$
is sc$_\triangleleft$-smooth. Similarly we can define sc-sections
and sc$^+$-sections $\Gamma({\mathcal
S}_{\triangleleft},{K}^{{\mathcal S}_{\triangleleft}}|{ O})$ and
$\Gamma^+({\mathcal S}_{\triangleleft}, { K}^{{\mathcal
S}_{\triangleleft}}|{O})$. In future, if irrelevant, we might
suppress the ${\mathcal S}_{\triangleleft}$ in the notation and
write for example $\Gamma({K}^{{\mathcal S}_{\triangleleft}}|{ O})$.
We also consider sc$_{\triangleleft}$-smooth vector bundle morphisms
which are those which are linear in the fiber.

{ Similarly as we introduced M-polyfolds and strong sc-vector
bundles, we can define M-polyfold bundles $b:Y\rightarrow X$. Let us
also remark that one can develop a good notion of connection for
$b$. These strong sc-connections have special properties reflecting
the fact that we have a grading of the fiber, say
$Y_{m,k}\rightarrow X_m$ with $0\leq k\leq m+1$ and the compact
inclusion from level $m+1$ to $m$. As a consequence covariant
derivatives of a section $f$ of $b$ with respect to different
connections in this class differ by a sc$^+$-operator. In particular
the difference is a linear sc-compact operator.} This will be
important for the finer aspects of the Fredholm theory as needed for
the "operation theory" in \cite{HWZ3}, f.e. orientation questions.

Coming back to our Morse-theory example we can define a M-polyfold
bundle $\overline{Y}$ over the M-polyfold $\overline{X}$. Its
elements are sequences $([h_1],...,[h_k])$ of equivalence classes of
$H^1$-sections $h_i$ along underlying curves $u_i$, where
$([u_1],..,[u_k])\in \overline{X}$. With the notion of sc-Fredholm
section, which we are going to introduce in the next section, it
will turn out the map
$$
f:\overline{X}\rightarrow \overline{Y}:([u_1],...,[u_k])\rightarrow
([\dot{u}_1-\Phi'(u_1)],..,[\dot{u}_k-\Phi'(u_k)])
$$
will be a sc-smooth Fredholm section. Moreover, for every section
$s\in\Gamma^+(\overline{Y})$ the section $f+s$ is always sc-Fredholm
for the bundle $\overline{Y}^1\rightarrow \overline{X}^1$ (Here all
indices are lifted by one.) Further, if $f$ is on every connected
component of $\overline{X}$ proper, the same will be true for $f+s$
if $s$ is small enough and has its support in a suitable open
neighborhood of $f^{-1}(0)$.

\subsection{Polyfold Groupoids and Polyfolds}
In dealing with Gromov-Witten theory or more generally with SFT we
need an orbifold version of the notion of M-polyfold. Orbifolds
always arise if we consider objects modulo some equivalence. Objects
with self-symmetries occur as singular points, i.e. true orbifold
points. There are different ways of defining orbifolds and similarly
(with some modifications) different ways of defining their
polyfold-generalization. The approach via groupoids seems
particularly useful, even from an analysis viewpoint. We refer the
reader to the excellent article by Moerdijk about a groupoid
approach to orbifolds, \cite{Moerdijk}, as well as the book
\cite{MoerMr}.

\subsubsection{Polyfold groupoids}
Recall that a groupoid $\mathfrak{G}$ is a small
category where every morphism is invertible. We shall write $G$ for
the objects and ${\bf G}$ for the morphism set. Given a groupoid we
have a certain number of obvious structure maps. There are the
source and target maps
$$
s,t:{\bf G}\rightarrow G
$$
which associate to a morphism its source or target. Since every
morphism is invertible we have the inversion map
$$
i:{\bf G}\rightarrow {\bf G}:\phi\rightarrow \phi^{-1}.
$$
In addition we have the unit map
$$
u:G\rightarrow {\bf G}:x\rightarrow 1_x.
$$
Finally we can build the fibered product ${\bf G}{_s\times_t}{\bf
G}$ consisting of all pairs of morphisms $(\phi,\psi)$ with
$s(\phi)=t(\psi)$. Then we can define the multiplication map
$$
m:{\bf G}{_s\times_t}{\bf G}\rightarrow {\bf
G}:(\phi,\psi)\rightarrow \phi\circ \psi.
$$
Now we are almost in the position to introduce the notion of a
polyfold groupoid. As a final preparation we need the notion of a
fred-submersion between two M-polyfolds $X$ and $Y$. If
$\mathcal{T}=(\pi,E,V)$ and $\mathcal{S}=(\rho,F,V)$ are splicings
with common parameter set $V$ we can build the Whitney sum
${\mathcal T}\oplus{\mathcal S}$ by defining
$$
{\mathcal T}\oplus{\mathcal S}=(\tau,E\oplus F,V),
$$
where
$$
\tau_v(e,f)=(\pi_v(e),\rho_v(f)).
$$
A particular situation arises if ${\mathcal S}=(Id,{\mathbb
R}^N,V)$. In that case we will simply write ${\mathcal T}\oplus
{\mathbb R}^N$ instead of ${\mathcal T}\oplus{\mathcal S}$.
\begin{definition}
A sc-smooth map $f:X\rightarrow Y$ between the M-polyfolds $X$ and
$Y$ is said to be a fred-submersion, if at every point $x_0\in X$
resp. $f(x_0)\in Y$ there exists a chart $(U,\varphi,{\mathcal
T}\oplus {\mathbb R}^N)$ resp. $(W,\psi,{\mathcal T})$ satisfying
$f(U)\subset W$ and
$$
\psi\circ f\circ \varphi^{-1}(v,e',e'') =(v,e').
$$
\end{definition}
Note the following easy consequence of the definition of a
fred-submer\-sion
\begin{proposition}
If $f:X\rightarrow Y$ is a fred-submersion between the M-polyfolds
$X$ and $Y$ then for every smooth $y\in Y$ the preimage $f^{-1}$
carries in a natural way the structure of a smooth
finite-dimensional manifold.
\end{proposition}
Given three M-polyfolds $X$, $X'$ and $Y$ and sc-smooth maps
$s:X\rightarrow Y$ and $t:X'\rightarrow Y$ we can build as a set the
fibered product $X{_s\times_t} X'$ by defining
$$
X{{_s}\times_t} X'=\{(x,x')\ |\ s(x)=t(x')\}.
$$
In certain situations this set carries in a natural way the
structure of a M-polyfold. We have
\begin{proposition}
If at least one of the maps $s$ or $t$ is a fred-submersion the
fibered product $X{{_s}\times_t} X'$ carries in a natural way the
structure of a M-polyfold. Further, if $s$ is a fred-submersion the
projection
$$
X{{_s}\times_t} X'\rightarrow X'
$$
is also a fred-submersion. If $t$ is a fred-submersion the same is
true for the projection
$$
X{{_s}\times_t} X'\rightarrow X.
$$
\end{proposition}
Now we can give the definition of a polyfold groupoid.
\begin{definition}
A polyfold groupoid is a groupoid $\mathfrak{X}$, together with a
M-polyfold structure for the set of objects $X$ and the set of
morphisms ${\bf X}$ so that the source and target maps $s$ and $t$
are surjective fred-submersions and all structure maps are
sc-smooth. We assume the induced topologies on $X$ and ${\bf X}$ to
be second countable and paracompact.
\end{definition}
Note that ${\bf X}{{_s}\times_t}{\bf X}$ is a M-polyfold since $s$
and $t$ are fred-submersions, so that it makes sense to talk about
the sc-smoothness of the multiplication map $m$.

Clearly, the notion of a polyfold groupoid, is a straight forward
modification of that of a Lie groupoid, where we have replaced the
notion of a finite-dimensional manifold by that of a M-polyfold and
the notion of submersion is modified by that of a fred-submersion,
see \cite{Moerdijk,MoerMr}. For the convenience of the reader let us
recall the  definition:
\begin{definition}
A Lie groupoid is a small category $\mathfrak{X}$, where the set of
objects $X$ and morphisms ${\bf X}$ is equipped with a smooth
manifold structure\footnote{For our purposes we can assume that the
manifolds are second countable. However, for certain applications
(not relevant for us) one should allow non-Hausdorff manifolds, see
\cite{MoerMr}.}, so that the source and target maps are surjective
submersions and all structure maps are smooth.
\end{definition}

 The orbit space
$|\mathfrak{X}|$ of a polyfold groupoid consists of the quotient
space $X/\sim$, where two points are identified if they are related
by a morphism. Observe that $|\mathfrak{X}|$ inherits a filtration
from $X$.  The maps between two polyfold groupoids are the sc-smooth
functors $F:\mathfrak{X}\rightarrow \mathfrak{Y}$. That means $F$
induces sc-smooth maps $X\rightarrow Y$ and ${\bf X}\rightarrow {\bf
Y}$. A sc-smooth functor $F$ induces a sc$^0$-map between the
orbit-spaces.

For the construction of polyfolds it will be important to introduce
the notion of a generalized map. The discussion is similar to that
in the Lie groupoid situation. For this let us first introduce the
notion of an equivalence.
\begin{definition}
Let ${\mathfrak{X}}$ and $\mathfrak{Y}$ be M-polyfold groupoids. A
sc-smooth functor $F:\mathfrak{X}\rightarrow\mathfrak{Y}$ is called
an equivalence provided the following holds:
\begin{itemize}
\item[1)] The map $t\pi_1:{\bf Y}{{_s}\times_F}X\rightarrow Y$ is a
surjective fred-submersion.
\item[2)] The square
\begin{eqnarray*}
\begin{CD}
{\bf Y} @>F>> {\bf X}\\
@VV{(s,t)}V     @VV{(s,t)}V\\
Y\times Y @>F\times F>>  X\times X
\end{CD}
\end{eqnarray*}
is a fibered product.
\end{itemize}
\end{definition}
An equivalence is usually not invertible. At this point we have a
category where the objects are M-polyfold groupoids with the
sc-smooth functors being the morphisms between them. Further we have
a distinguished family of special morphisms, namely the
equivalences. There is now a very particular, purely
category-theoretic procedure for inverting a distinguished class of
arrows in a category, while at the same time keeping the objects and
only minimally change (given the fact that we must invert a certain
number of given arrows) the morphisms. The general procedure is
described in \cite{GZ}. Here we will describe the procedure for our
special situation. We need a certain amount of preparation.
\begin{definition}
Assume that $F,G:\mathfrak{A}\rightarrow\mathfrak{B}$ are sc-smooth
functors. They are called equivalent if there exists a sc-smooth map
$$
\tau:A\rightarrow {\bf B}
$$
associating to every object $x\in A$ a morphism
$$
\tau(x):F(x)\rightarrow G(x)
$$
so that for every $h:x\rightarrow x'$ we obtain the commutative
diagram
\begin{eqnarray*}
\begin{CD}
F(x) @>\tau(x)>> G(x)\\
@VVF(h)V     @VVG(h)V\\
F(x') @>\tau(x')>> G(x').
\end{CD}
\end{eqnarray*}
The map $\tau$ is called a natural transformation.
\end{definition}
In order to define generalized maps start with a diagram
$$
\mathfrak{X}\xleftarrow{F}\mathfrak{A}\xrightarrow{\phi}\mathfrak{Y},
$$
where $F$ is an equivalence and $\phi$ a sc-smooth functor. Let us
call $\mathfrak{X}$ the domain and $\mathfrak{Y}$ the codomain (of
the diagram). Consider a second such diagram
$$
\mathfrak{X}\xleftarrow{F'}\mathfrak{B}\xrightarrow{\phi'}\mathfrak{Y},
$$
with identical domain and codomain. We call it a refinement of the
first if there exists a sc-smooth functor
$H:\mathfrak{B}\rightarrow\mathfrak{A}$ so that $F\circ H$ and $F'$
are naturally equivalent as well as $\phi\circ H$ and $\phi$.
Finally we say that two diagrams of maps, say
$$
\mathfrak{X}\xleftarrow{F}\mathfrak{A}\xrightarrow{\phi}\mathfrak{Y}\
\hbox{and}\
\mathfrak{X}\xleftarrow{F'}\mathfrak{A'}\xrightarrow{\phi'}\mathfrak{Y}
$$
are equivalent if they have a common refinement. It takes a certain
amount of work to show that this indeed defines an equivalence
relation. We associate to a smooth functor
$\phi:\mathfrak{X}\rightarrow\mathfrak{Y}$  the equivalence class of
$$
\mathfrak{X}\xleftarrow{Id}\mathfrak{X}\xrightarrow{\phi}\mathfrak{Y}.
$$
Let us denote this equivalence class by $[\phi]$. Similarly we
denote for an equivalence $F:\mathfrak{X}\rightarrow \mathfrak{Y}$
the class of the diagram
$$
\mathfrak{Y}\xleftarrow{F}\mathfrak{X}\xrightarrow{Id}\mathfrak{X}
$$
by $[F]^{-1}$. Then one verifies immediately that this is the
inverse of $[F]$.

\subsubsection{Polyfolds}
Now we are in the position to define polyfolds. These type of spaces
suffice to carry out all the functional analytic constructions in
Gromov-Witten, Floer-Theory, and SFT.

\begin{definition} Consider a polyfold groupoid $\mathfrak{X}$.
\begin{itemize}
\item[1)] We say that $\mathfrak{X}$ is \'etale provided the source
and target maps are local sc-diffeomorphisms.
\item[2)] We say that $\mathfrak{X}$ is proper if for every $x\in X$
there exists an open neighborhood $U(x)$ so that the map
$$
t:s^{-1}(\overline{U(x)}\rightarrow X
$$
is proper.
\item[3)] A polyfold groupoid which is
\'etale and proper is called an ep-polyfold groupoid.
\end{itemize}
\end{definition}
Note that we could interchange the role of $s$ and $t$ in the
definition of proper defining the same property.
\begin{definition}
Let $Z$ be a second countable paracompact topological space. An
ep-polyfold structure on $Z$ is given by a pair
$(\mathfrak{X},\beta)$, where $\mathfrak{X}$ is a ep-polyfold
groupoid and $\beta:|\mathfrak{X}|\rightarrow Z$ a homeomorphism. We
say that two polyfold structures $(\mathfrak{X},\beta)$ and
$(\mathfrak{X}',\beta')$ are equivalent if there exist equivalences
$$
F:\mathfrak{X}''\rightarrow\mathfrak{X}\ \hbox{and}\
\mathfrak{X}''\rightarrow \mathfrak{X}'
$$
so that
$$
\beta\circ |F|=\beta'\circ |F'|.
$$
\end{definition}
Let us observe that $(\mathfrak{X}'',\beta\circ F)$ is also a
polyfold structure equivalent to the two other ones. Finally we can
introduce the notion of a polyfold
\begin{definition}
A second countable paracompact topological space $Z$ equipped with
an equivalence class of polyfold structures is called a polyfold.
\end{definition}
For a polyfold $Z$ we can define a degeneracy function
$d:Z\rightarrow {\mathbb N}$ by
$$
d(z)=d_X(x)
$$
where $(\mathfrak{X},\beta)$ is a defining polyfold structure and
$\beta(x)=z$. Here $d_X$ is the degeneracy map for $X$. This is
well-defined, i.e. independent of the representative of the polyfold
structure. We leave it to the reader to verify that it makes sense
to talk about a face-structured polyfold.

\subsubsection{Polyfold bundles}
We have already defined strong bundles over M-polyfolds. In a first
step we introduce strong bundles over a polyfold groupoid. Let us
start with a given polyfold groupoid $\mathfrak{X}$ and a strong
M-polyfold bundle $\tau:E\rightarrow X$ over $X$. Using the fact
that the source map $s:{\bf X}\rightarrow X$ is a fred-submersion we
can build the pull-back bundle via $s$ over ${\bf X}$. This
pull-back bundle is, of course,
$$
{\bf X}{_s\times_\tau}E\rightarrow {\bf X}.
$$
This is a strong bundle over ${\bf X}$. Assume we are given a strong
bundle map $\mu:{\bf X}{_s\times_\tau}E\rightarrow E$ covering $t$.
To be more precise $\mu$ gives the following commutative diagram
\begin{eqnarray*}
\begin{CD}
{\bf X}{_s\times_\tau}E @>\mu>> E\\
@VV\pi_1 V     @VV\tau V\\
{\bf X} @>t>> X.
\end{CD}
\end{eqnarray*}
We require now $\mu$ to be compatible with the morphisms in ${\bf
X}$. More precisely we require with the abbreviation $g\cdot
e=\mu(g,e)$ that
\begin{itemize}\label{propxx1}
\item The identity $1_x\cdot e=e$ holds.
\item Moreover $(g\circ h)\cdot e = g\cdot(h\cdot e)$.
\end{itemize}
The following definition is useful:
\begin{definition}
Let $\mathfrak{X}$ be a polyfold groupoid. A strong linear
$\mathfrak{X}$-space is given by a pair  $(E,\mu)$, where
$\tau:E\rightarrow X$ is a strong bundle over $X$ and $\mu: {\bf
X}{_s\times_\tau}E \rightarrow E$ a strong vector bundle isomorphism
satisfying the properties in (\ref{propxx1}).
\end{definition}
Given a strong linear $\mathfrak{X}$-space $(E,\mu)$ we can build a
category $\mathfrak{E}$, with objects $E$ and morphism set ${\bf E}$
defined by
$$
{\bf E}:={\bf X}{_s\times_\tau}E.
$$
We define the source map $s$ by $s(g,e)=e$ and the target map by
$t(g,e)=g\cdot e$. This defines a polyfold groupoid $\mathfrak{E}$.
Note however that we have more structure since $E\rightarrow X$ for
example is a strong bundle. The projection map $\tau:E\rightarrow X$
extends to a sc-smooth functor
$$
\tau:\mathfrak{E}\rightarrow \mathfrak{X}.
$$
We might view the latter diagram as the strong linear
$\mathfrak{X}$-space, keeping the extra structure of $\mathfrak{E}$
in mind.

Next we introduce polyfold bundles. Assume we start with a
surjective continuous map $\tau:L\rightarrow  Z$ between second
countable paracompact spaces. We assume in addition that the fibers
are equipped with Banach space structures. Consider a strong bundle
$\pi:\mathfrak{E}\rightarrow\mathfrak{X}$ over the ep-polyfold
groupoid $\mathfrak{X}$ and assume we are given homeomorphisms
$$
\beta:|\mathfrak{E}|\rightarrow L\ \hbox{and}\ \
\beta_0:|\mathfrak{X}|\rightarrow Z,
$$
so that
$$
\tau\circ\beta=\beta_0\circ|\pi|.
$$
We assume that $\beta$ is linear on the fibers. We can define
equivalence classes of such objects
$(\pi:\mathfrak{E}\rightarrow\mathfrak{X},\beta,\beta_0)$ as before
by refinements. Of course linearity should be preserved in the
fibers (the precise details we leave to the reader).
\begin{definition}
Assume that we are given two second countable paracompact spaces
together with a surjective continuous map $\tau:L\rightarrow Z$. In
addition every fiber is equipped with a Banach space structure whose
topology coincides with the given topology. A strong polyfold bundle
structure for $\tau$ is given by an equivalence class of triples
$(\pi:\mathfrak{E}\rightarrow\mathfrak{X},\beta,\beta_0)$ as
described above.
\end{definition}

A sc-smooth section $f$ of $\tau:L\rightarrow Z$ would we be
represented by a sc-smooth section (functor) $F$ of some
$\pi:\mathfrak{E}\rightarrow\mathfrak{X}$, where the latter occurs
in a triple defining the strong polybundle structure. We would like
to point out that these abstract bundles together with the Fredholm
theory developed in the following section are sufficient as a
functional analytic framework for Gromov-Witten theory. Our Fredholm
theory will allow for an abstract multi-valued perturbation theory
so that the solutions spaces in Gromov-Witten theory are smooth
compact branched manifolds. McDuff has developed a convenient
frame-work, \cite{McDuff}, in which integration, differential forms
etc. make sense so that many formulas in Gromov-Witten theory can be
in fact obtained through integration of suitable quantities over the
moduli spaces. So far the smoothness theory (in the symplectic case)
for the moduli spaces was not developed  to the extend necessary.
\subsection{Comments}
Let us finish the section with some remarks about possible
variations of our previous definitions. For example one could remove
the condition that the embedding $E_n\rightarrow E_m$ is compact for
$n>m$. In that case the sequence $E_m=E$ would be allowed. However,
in order to obtain the chain rule the notion of a sc-smooth map has
to be modified as well. In fact one has to replace the requirement
that
$$
U_{m+1}\oplus E_m\rightarrow F_m:(x,h)\rightarrow Df(x)h
$$
is continuous by the requirement that
$$
U_{m+1}\rightarrow L(E_m,F_m):x\rightarrow Df(x)
$$
is continuous. In the case of the constant sequence $E_m=E$ we
recover the standard notion of Frechet differentiability in a Banach
space. As it turns out this modified concept and the associated
theory parallel to the previous discussion is neither applicable to
Gromov-Witten theory, Floer Theory and SFT. Observe that our notion
of smoothness and the above modified notion only coincide in the
finite-dimensional case, but are otherwise strikingly different.
Hence the notion of differentiability which works in our
applications of interest generalizes the finite-dimensional notion,
but does not contain Frechet-differentiability in
infinite-dimensional Banach spaces (which is the commonly used
generalization) as a special case.

One can also generalize the notion of splicing allowing for an
infinite-dimensional set of splicing parameters, but the author is
not currently aware of any good application.

The notion of ep-polyfold can be generalized by not requiring the
\'etale condition. In this case one could study objects with a
compact Lie group as a symmetry group (isotropy group). This would
be a necessary generalization to deal with Yang-Mills type problems.

\section{Generalized Fredholm Theory }
In this section we describe the Fredholm theory.  The article
\cite{BSV} gives an overview over the classical Fredholm theory. The
constructions and ideas described in \cite{BSV} can be carried out
within our generalized Fredholm context. In addition many more
constructions and concepts are possible and lead to a theory
applicable to a much wider range of problems. Moreover, the
"Fredholm Theory with Operations" which we describe in some simple
cases later, gives a broad range of new structures on Fredholm
problems. In principle one could have formulated such an abstract
theory also within the  classical Fredholm theory. Unfortunately,
these structures have only been observed as consequences of the lack
of compactness and clever compactifications of the moduli spaces,
i.e. in situations where the classical theory is not applicable.

We forget sc-structures for the moment and have a look at the
classical Fredholm-situation. Assume that $f:U\rightarrow F$ is a
smooth map, which is defined on the open neighborhood $U\subset E$
of $0$ and takes values in $F$. Moreover $f(0)=0$. Here $E$ and $F$
are Banach spaces. We assume that $f'(0)$ is a linear Fredholm
operator. Then there exist topological splittings of the domain
$E=K\oplus X$ and the target $F=C\oplus Y$ with the property that
$$
f'(0):X\rightarrow Y
$$
is a linear isomorphism. Define an isomorphism
$$
\sigma:F=C\oplus Y\rightarrow C\oplus X:(c,y)\rightarrow
(c,f'(0)^{-1}(y)).
$$
Then the composition
$$
K\oplus X\rightarrow C\oplus X:(k,x)\rightarrow \sigma\circ
f'(0)(k+x)
$$
has the form $(k,x)\rightarrow (0,x)$. If we consider now
$\sigma\circ f(k,x)$ and project in the target onto $X$ along $C$ we
obtain a map of the form
$$
(k,x)\rightarrow x-B(k,x)
$$
where for $k$ small $x\rightarrow B(k,x)$ is a contraction. With
other words, by appropriately taking a suitable coordinate
representation of $f$, a suitable splitting of the domain and  a
projection onto a finite codimension subspace this new map looks
like a parameterized contractive perturbation of the identity. If
$f(x_0)=y_0 \neq 0$ and we want to understand the behavior near
$x_0$ then we can look at $x\rightarrow g(x):=f(x+x_0)-f(x_0)$ which
brings us back into the first case. It is not difficult to verify
that the contraction "normal form" together with a smoothness
requirement is equivalent to saying that the original $f$ is
Fredholm. This equivalent formulation fits well into our
sc-framework as can be seen in the following.

\subsection{Contraction Germs}

Now we introduce the relevant result for our situation.  If we
talk about a sc-germ $f:\mathfrak{O}({E},0)\rightarrow ({ F},0)$
we mean that for every level $m$ we have an open neighborhood
around $0$ (for the topology on the $m$-level) on which $f$ is
defined and maps it into $m$-level. Below ${\mathbb N}$ denotes
the non-negative integers.

\begin{definition}
Let $f:\mathfrak{O}(V\oplus {E},0)\rightarrow ({ E},0)$ be a
sc$^0$-germ with $f(0,0)=0$, where $V$ is an open subset of a
partial cone in some finite-dimensional vector space and ${ E}$ is
a sc-Banach space. We call $f$ a sc-contraction germ if $f$ can be
written in the form
$$
f(v,u)=u-B(v,u)
$$
so that the following holds. For every level $m\in {\mathbb N}$,
and a suitable $\Theta_m\in (0,1)$ we have an estimate
$$
\parallel B(v,u)-B(v,u')\leq \Theta_m\cdot\parallel u-u'\parallel_m
$$
provided $v,u,u'$ are close enough to $0$ (the notion of close
depending on $m$ and $\Theta_m$)
\end{definition}

Banach's fixed point theorem applied to every level gives a
sc$^0$-germ $\delta:\mathfrak{O}(V,0)\rightarrow
\mathfrak{O}(E,0)$ so that its graph $gr(\delta)$ satisfies
$$
f\circ gr(\delta)=0.
$$
The main result is the following "Germ-Implicit Function Theorem":
\begin{theorem}[Germ-Implicit Function Theorem]\ \
Let $f:\mathfrak{O}(V\oplus E,0)\rightarrow \mathfrak{O}(E,0)$ be
sc-smooth and a sc$^0$-contraction germ. Then the solutions germ
$\delta$ of $f\circ gr(\delta)=0$ is sc-smooth.
\end{theorem}
To be more precise, the conclusion is that for every $m$ and $k$ and
$|v|$ small enough the map $v\rightarrow \delta(v)$ goes into the
$m$-level and is $C^k$. In particular $\delta$ is $C^{\infty}$ at
the point $0$. As it turns out, describing more globally the
solution set of a problem of the form $f=0$, all these local
solution germs fit smoothly together to give the solution set a
smooth structure. Via this observation the above theorem will be one
of the key building blocks for all versions of the implicit function
theorem, as well as transversality theory.

\subsection{Fillers}
We would like to use the above discussion to generalize the idea of
a Fredholm section in our polyfold set-up. The study above only
centers at a particular local situation which takes place in open
subsets of Banach spaces (with a sc-structure). Clearly we need to
explain how splicings come in if we want to develop a theory in
polyfolds. Of course, one might expect that locally varying
dimensions of the ambient spaces lead to a cumbersome definition of
what a Fredholm operator in such a context would really mean.
Surprisingly the idea of a "Filler" which is a rather simple object
makes it a non-issue.

We start with the local set-up. Assume that we have a strong local
sc-vector bundle
$$
U\triangleleft F\rightarrow U,
$$
where $U$ is open in the sc-Banach space $E$. Since we will only be
interested in the neighborhoods of smooth points we may without loss
of generality assume that $0\in U$. Also, being only interested in a
neighborhood of $0$ we will write $\mathfrak{O}(U\triangleleft F,0)$
to emphasize that we consider the germ of a bundle. Given a germ of
a section $[f,0]$ we denote by $\bar{f}$ the principal part of $f$.
Then
$$
\bar{f}:\mathfrak{O}(U,0)\rightarrow \mathfrak{O}(F,\bar{f}(0)).
$$
Linearize $\bar{f}$ at $0$, say $\bar{f}'(0):E\rightarrow F$. We say
that $[f,0]$ is linearized Fredholm if $\bar{f}'(0)$ is sc-Fredholm.
Observe, that the linearization of a section $f$ of some abstract
vector bundle at some point $q$ with $f(q)\neq 0$ is not an
intrinsic object, whereas it is at a zero. It depends on the choice
of local coordinates.

In our case, however, due to the notion of strong bundle the
property of being linearized Fredholm is in fact intrinsic in the
following sense. If $\Psi:\mathfrak{O}(U\triangleleft
F,0)\rightarrow \mathfrak{O}(V\triangleleft G,0)$ is a germ of a
sc-vector bundle isomorphism, then the push-forward germ
$\Psi_{\ast}([f,0])$ is linearized Fredholm if and only this is true
for $[f,q]$. This is a deeper consequence of the property of a
strong(!) sc-vector bundle. The reader might verify, that the
linearizations taken of two different local coordinate
representations (using strong bundle coordinates), say $L_1$ and
$L_2$ are related by
$$
L_1 = AL_2B +K,
$$
where $A$ and $B$ are sc-isomorphisms and $K$ is a sc$^+$-operator.
In particular $K$ is level-wise a compact perturbation. Hence, given
a section of a strong sc-vector bundle $b:Y\rightarrow X$, saying
that for a $q\in X^{\infty}$, i.e. a smooth $q$, the germ $[f,q]$ is
linearized Fredholm has an intrinsic meaning (of course the actual
linearization depends on the choice of local coordinates). Up to
this point we have discussed (germs of) sections of a strong bundle
over a sc-manifold.

Let us assume next that we have a germ of sc-smooth section $[f,0]$
of a local M-polyfold bundle. With other words for $(v,e)$ in the
splicing core $K^{{\mathcal S}_0}$ near $0$ the map $f$ takes the
form
$$
f(v,e)=((v,e),\bar{f}(v,e))
$$
with image in $K^{{\mathcal S}_0}\triangleleft_V K^{{\mathcal
S}_1}$. Let us abbreviate $K_i=K^{{\mathcal S}_i}$. If we fix $v$
the space $K_{0,v}=\{e\ |\ \pi_v(e)=e\}$ is a sc-Banach space and
similarly $K_{1,v}$. The principal part $\bar{f}$ has the property
that $e\rightarrow \bar{f}(v,e)$ maps $K_{0,v}$ to $K_{1,v}$. We say
that $[f,0]$ is linearized Fredholm provided the derivative of map
$e\rightarrow \bar{f}(0,e)$ at $0$ is sc-Fredholm. It turns out that
this is again an intrinsic definition invariant under changes of
coordinates.

If we have a splicing ${\mathcal S}=(\pi,E,V)$ we also have the
complementary splicing ${\mathcal S}^c=(I-\pi,E,V)$. Clearly
$V\oplus E$ can be viewed as the fiber sum over $V$
$$
V\oplus E= K^{\mathcal S}\oplus_V K^{{\mathcal S}^c}.
$$
In particular a point $(v,e)$ in $V\oplus E$ can be written as
$$
(v,e) = (v,u_v+u_v^c),
$$
where $v\in V$, $u_v\in K^{\mathcal S}_v$ and $u_v^c\in K^{{\mathcal
S}^c}_v$. A section germ $[f,0]$ for $\mathfrak{O}(K^{{\mathcal
S}_0}\triangleleft_V K^{{\mathcal S}_1},0)$ is called fillable if
there exists a section germ $[\hat{f},0]$ of $\mathfrak{O}((V\oplus
E)\triangleleft F,0)$ having the form
$$
\hat{f}(v,u_v+u^c_v)=((v,u_v+u^c_v),\bar{f}(v,u_v)+\bar{f}^c(v,u_v,u^c_v)),
$$
where $\bar{f}^c$ is defined on an open neighborhood of $(0,0)$ in
$V\oplus E$ mapping its points to $K^c_1$ in such a way that
$u^c_v\rightarrow \bar{f}^c(v,u_v,u^c_v)$ is a linear sc-isomorphism
$K^c_{0,v}\rightarrow K^c_{1,v}$.

What is the significance of a filler? We would like to study the
section $f$, meaning that we are interested in the solution set of
$f=0$. Here, for $(v,u_v)\in K^{{\mathcal S}_0}_v$ we have
$(v,\bar{f}(v,u_v))\in K^{{\mathcal S}_1}_v$. Consider now the
filled section $\hat{f}$. If $\hat{f}(v,u)=0$ we conclude that
$$
\bar{f}^c(v,u_v,u_v^c)=0.
$$
By the properties of the filler this means that $u_v^c=0$. Hence we
conclude that $(v,u)=(v,u_v)\in K^{{\mathcal S}_0}_v$ and
$$
\bar{f}(v,u)=0.
$$
With other words, the modification by a filler does not change the
solution set. Hence, locally the study of $\bar{f}$, which is
defined on a perhaps very bad space with varying dimensions is
equivalent to the study the section $\hat{f}$ on open sets of Banach
spaces with sc-structure. The nice fact is in applications, that
there are usually obvious choices for fillers. For example in
Morse-theory every critical point $b$ has an associated filler
$$
h\rightarrow \dot{h}-\Phi''(b)h,
$$
where $h:{\mathbb R}\rightarrow T_bM$ belongs to a suitable
sc-Hilbert space of functions and $\Phi''(b)$ is the Hessian. In
Gromov-Witten theory the situation is slightly more complicated. The
fillers are associated to the images of the nodal points which can
be any point on the symplectic manifold $W$. Then the filler
associated to $w\in W$ is
$$
h\rightarrow h_s+J(w)h_t,
$$
i.e. the linear Cauchy-Riemann operator acting on maps
$$
h:{\mathbb R}\times S^1\rightarrow T_wW,
$$
where $h$ takes antipodal values at $\pm\infty$:
$$
h(-\infty)+h(+\infty)=0
$$
and belongs to a certain Sobolev class. In SFT the asymptotic
periodic orbits behave like  one-dimensional Morse-Bott manifolds of
critical points and the fillers are one-dimensional families of
(linear) perturbed Cauchy-Riemann type problems associated to
periodic orbits. We refer the reader to a comprehensive discussion
of fillers in \cite{HWZ1} including a complete discussion of the
Morse-theory case. The fillers for SFT are constructed in full
detail in \cite{HWZ2}.

Let us conclude, that the bottom line of the "filler discussion" is
that we can reduce the local study of sections of strong M-polyfold
bundles to the study of sections of strong local bundles. In this
context the previous results on "contraction germs" allows to derive
suitable implicit function theorems.

\subsection{Fredholm Operators}
{ Let $b:Y\rightarrow X$ be a smooth M-polyfold bundle and $f$ a
section. We will define the notion of a  Fredholm section. We give a
less condensed form (than is possible) since it is more instructive.
First of all $f$ is sc-smooth and regularising, i.e. if $f(x)\in
Y_{m,m+1}$ then $x\in X_{m+1}$. Of course, if $f$ is regularizing
and $f(x)=0$ we can conclude that $x\in X_{\infty}$. Let us also
note the following important fact for the perturbation theory using
sc$^+$-sections, $s\in\Gamma^+(b)$. If $f(x)+s(x)\in Y_{m,m+1}$ for
$s\in\Gamma^+(b)$, then $x\in X_{m}$ and $s(x)\in Y_{m,m+1}$. The
latter is true by definition of a sc$^+$-section. Consequently
$f(x)\in Y_{m,m+1}$ implying that $x\in X_{m+1}$. With other words
$f+s$ is also regularizing.

Secondly, for every smooth $q\in X$, i.e. $q\in X_{\infty}$, there
is a local strong M-polyfold bundle trivialization mapping the germ
$[f,q]$ to a fillable $[f_1,0]$. Further there exists a filled
section $[\hat{f},0]$ which is a section germ of
$\mathfrak{O}((V\oplus E)\triangleleft F,0)$. Thirdly there is a
sc-Banach space $W$, a finite dimensional vector space $R$, another
finite-dimensional vector space $Q$ and an open subset $B$ of some
partial cone in $Q$ and a local strong sc$_\triangleleft$-vector
bundle isomorphism
$$
\Psi:\mathfrak{O}((V\oplus E)\triangleleft F,0)\rightarrow
\mathfrak{O}((B\oplus W)\triangleleft (R\oplus W),0)
$$
so that the push-forward $[g,0]$ of $[\hat{f},0]$ has the
following property. If $P:R\oplus W\rightarrow W$ is the
projection and $\bar{g}$ the principal part, then
$$
(b,w)\rightarrow P\bar{g}(b,w)-P\bar{g}(0,0)
$$
is a sc$^0$-contraction germ. In other words:

\noindent {\it "A section of a M-polyfold bundle is Fredholm
provided it is regularizing and in suitable local coordinates it
admits a filler, so that the filled section gives under another
coordinate change and suitable splittings and projection a
contraction germ."}

Let us denote by $\hbox{Fred}(b)$ the Fredholm sections of $b$.} The
definition of a Fredholm section is very general. It looks not very
practical at first sight, but at least as applications are
concerned, it indeed is and the method for showing that the
nonlinear elliptic pde's in GW, FT, CH and SFT are Fredholm in our
generalized sense are almost identical. Let us elaborate somewhat
about this point. The advantage of the classical implicit function
theorem is, of course, the fact that we can conclude something about
the local properties of the solution set by knowing something about
the linearization at a single point. In our applications it is,
however, not applicable. Nevertheless in practice our problems will
in general only have a finite-dimensional set of bad parameters
which will not enter smoothly (but they enter in a sc-smooth way).
With respect to the remaining variables we will have smoothness and
linearizations (in the classical way on every level). A certain
uniformity of behavior of these linearizations with respect to the
bad parameters allows to show the "contraction normal form". Hence,
in applications, the analysis needed is concerned with the standard
implicit function theorem applied to continuous (finite-dimensional)
families of maps and the normal form is a consequence of certain
uniformity of the estimates. This is reminiscent of the uniformity
of estimates in the gluing constructions occurring in FT or GW, see
f.e. \cite{MS}. We also refer the reader to \cite{HWZ1}, where the
Fredholm property is shown in the context of Morse-theory,
validating our previous remarks.

There are quite a number of consequences of this definition. For
example if $f$ is a Fredholm section of $b:Y\rightarrow X$, and $s$
is a sc$^+$-section then $f+s$ is a Fredholm section of
$b^1:Y^1\rightarrow X^1$. Hence
\begin{proposition}
There is a well-defined map
$$
\hbox{Fred}(b)\times \Gamma^+(b)\rightarrow
\hbox{Fred}(b^1):(f,s)\rightarrow f+s.
$$
\end{proposition}

An important consequence of the germ-implicit function theorem and
the definition of Fredholm section is the following implicit
function theorem
\begin{theorem}
Assume that $b:Y\rightarrow X$ is a smooth M-polyfold bundle and $f$
a Fredholm section. Suppose further that $\partial X=\emptyset$. If
$f(q)=0$ and the linearization $f'(q):T_qX\rightarrow Y_q$ is onto
then the solution set of $f(x)=0$ is near $q$ a smooth manifold (in
a natural way).
\end{theorem}

\subsection{Comments}
There are also results concerning the case where $X$ has a
boundary with corners, i.e. $\partial X\neq\emptyset$. We refer the
reader for more results to \cite{HWZ1}.

 One can show that if $f$ is Fredholm there are many small
perturbations by sc$^+$-sections $s$ so that for every solution of
$f(x)+s(x)=0$ the linearisation is onto, i.e. $(f+s)^{-1}(0)$ is a
smooth manifold. Moreover, if $f$ is proper then $f+s$ will be
proper if $s$ is small enough in a suitable sense. In addition
generic perturbations put the solution set into general position to
$\partial X$, see \cite{HWZ1}.

In the case of polyfold bundles where the section is represented by
a smooth section functor $F$ of
$\mathfrak{E}\rightarrow\mathfrak{X}$ the compatibility with the
morphisms asks for a perturbation theory respecting the morphisms.
As is well-known this required compatibility obstructs
transversality in general, if we are allowed only single-valued
perturbation. This changes if multi-valued perturbations are
allowed. In this case we obtain for generic multi-valued
sc$^+$-perturbation as a solution set a smooth branched weighted
manifold with boundary with corners (in general position to the
boundary) as defined by D. McDuff, see \cite{McDuff}. For details of
the Fredholm theory in this context we refer the reader to
\cite{HWZ1}.

\section{Operations} In this section we describe the important
theory of operations, which allows to capture the "algebra"
underlying the structure of having infinitely many interacting
Fredholm operators. We will describe only a somewhat simplified
version in order to expose the ideas. We refer the reader to
\cite{HWZ3} for extensions which are in fact necessary for the
applications we have in mind. We begin with a very useful
rudimentary algebraic structure.

\subsection{Degeneration Structures} A set with relators is a pair
$(S,R)$ with a set $S$ and a subset $R$ of $S\times S\times S$. We
write $(A,B;C)$ for an element in $R$ and call $A$ the left-source,
$B$ the right-source and $C$ the target. Given $(S,R)$ consider
$(k+1)$-tuples $z_k=(A_0,..,A_k)$ with $A_i\in S$. A $1$-step
degeneration is a diagram
$$
z_k\leftarrow z_{k-1}
$$
where $z_k$ is obtained from $z_{k-1}$ by replacing an occurring
element $C$ by two elements $(A,B)$ if $(A,B;C)\in R$. A "Short
Degeneration Sequence" has the form $z_0\rightarrow z_1\rightarrow
z_2$.
\begin{definition}
A degeneration structure $\mathfrak{S}=(S,R)$ consists of a set with
relators $(S,R)$   so that the following properties hold:
\begin{itemize}
\item[1)] (Degeneration Finiteness) Given $Z\in S$ the number of
degeneration sequences starting at $(Z)$ is finite.

\item[2)] (Associativity) The set of short degeneration sequences
having prescribed target and source  is either empty or consists of
precisely two elements. In the latter case these have the form
$$
(Z)\rightarrow(A,B)\rightarrow(A,I,E)\ \hbox{and}\ \ (Z)\rightarrow
(A^\ast,E)\rightarrow (A,I,E).
$$
\item[3)] (Minimality) If $(A,B;C)$ and $(A',B';C)$ belong to $R$
and either $A=A'$ and $A$ is not decomposable or if $B=B'$ and $B$
is not decomposable then $(A,B)=(A',B')$.
\end{itemize}
\end{definition}
Here $A$ is called decomposable if there exist $X,Y\in S$ with
$(X,Y;A)\in R$. Observe that for a degeneration structure
$(A,B;C)\in R$ implies $A\neq C$ and $B\neq C$. This is a
consequence of the finiteness axiom. It is important to note that
\begin{proposition}
If there exists a degeneration sequence
$$
z_n\leftarrow...\leftarrow z_0=(Z)
$$
Then there exist exactly $n!$ degeneration sequences connecting
$z_0$ with $z_n$.
\end{proposition}
Degeneration structures will be used to organize large families of
interacting Fredholm problems. Here is an example how they might
occur. If $f$ is a Fredholm section of the polyfold bundle
$Y\rightarrow X$ we might consider the set $S$ of connected
components of $X$, i.e. $S=\pi_0(X)$.  All symplectic problems we
have mentioned have the following structure as Fredholm problems.
First of all the underlying strong bundle $\tau:Y\rightarrow X$ is
defined over a face-structured polyfold. The compactified moduli
space is defined by $f=0$ for a  suitable Fredholm section $f$ of
$\tau$. In general $X$ has infinitely components and the moduli
space in every component is compact, whereas the overall solution
space is not compact. Here is now a rough description of a common
feature of FT, CH, SFT and in fact Morse Theory viewed in an
appropriate way. If we are given a point $z$ which solves the
Fredholm problem $f(z)=0$ and $z$ belongs to a face of $Z$ then $z$
can be viewed as the "product" $z'\circ z''$ of two solutions $z'$
and $z''$ in different components of $X$ (For example a boundary
face of the moduli space of gradient lines from $a$ to $c$ consist
of the broken trajectories $(z',z'')$ factoring over the same
intermediate critical point $b$. With other words  $z'$ is a perhaps
broken gradient line connecting $a$ with an intermediate critical
point $b$ and $z'$ connects $b$ with $c$. Hence we may view $z$ as
the product $z'\circ z''$.).  Coming back to the abstract situation,
if one of them, say $z'$ is again a boundary point, the point $z'$
is again a product and so on. Here, of course, product refers to
some kind of composition law how to build out of two solutions a new
one. In general the way to write $z$ as a product is not unique and
the non-uniqueness depends on the degeneracy $d(z)$ (For example a
two-broken gradient line can be viewed in precisely two different
ways as product. Moreover, we can view it as triple product and, of
course, this example already shows an associativity property which
should be required in general.).  Moreover, in general, given two
solutions there might in fact be finitely many different recipes to
produce new solutions, i.e. the product is in fact multi-valued.
This "product structure" satisfies some basic axioms which are
common to all examples and we in fact do not need to know more  in
order to develop our general theory. We will call the rule or method
of producing out of two solutions a new one an operation. This will
be discussed in more detail later. The operation will make sense not
only for solutions but for elements in the ambient space $X$, or the
bundle, as well. In short given a point $a$ in the component $A$ and
a point $b$ in the component $B$ we can produce a point $z$ in some
component $Z$ provided the components $A$, $B$ and $Z$ are related,
i.e.
$$
(A,B;Z) \ \ \hbox{is a relator}.
$$
Further $z$ will belong to a face of $Z$ which we will denote by
$[A,B;Z]$. We will make this more precise in the next subsection.
Before that let us describe a little bit more the landscape of
degeneration structure and related concepts. There are in addition
to degeneration structures notions like degeneration modules.
Degenerations modules compare to degenerations structures as modules
compare to rings. Degeneration modules occur when organizing
homotopies of inter-depending families of Fredholm operators. For
example organizing holomorphic curves in symplectic cobordisms
between contact manifolds we have a positive and negative boundary
component. Degeneration structures help to organize the Fredholm
problems associated to the boundary components and a degeneration
module (in fact a bi-module over the other two structures) organizes
the Fredholm operators associated to the cobordisms.  We refer the
reader to \cite{HWZ3} for more details, and \cite{EGH} for some
inspiration where these concepts might show up (look out for
bi-modules in the usual sense).

Let us mention one algebraic aspect of degeneration structures.
Assume that $\Lambda$ is a ring. Consider the group $C(S,\Lambda)$
of maps $S\rightarrow \Lambda$. Define the convolution
$\alpha\ast\beta$ by
$$
(\alpha\ast\beta)(C)=\sum_{(A,B;C)\in R} \alpha(A)\beta(B)
$$
Then the properties of a degeneration structure imply that this is
well-defined and the convolution is associative. We will explain
this later on our Morse-theory example. Applying a similar procedure
to degeneration modules we will obtain a bi-module (in the usual
sense) over the previously constructed rings (for the right and left
Fredholm problem). See the introduction to Section \ref{section5}
for some suggestive formulas. Now we are in the position to define
operations.

\subsection{Operations on M-Polyfolds} Let $\pi:Y\rightarrow X$ be
a M-polyfold bundle over a face-structured M-polyfold and $(S,R)$ a
degeneration structure. Let us begin with a process which we might
call indexing. Denote the set of connected components of $X$ by $S$
and assume we are given a degeneration structure $(S,R)$ where the
set of relators is in 1-1 correspondence with the set of faces of
$X$. More precisely, for every connected component $Z\in S$ the
faces of $Z$ are in 1-1 correspondence with the subset $R_Z$ of $R$
consisting of all relators of the form $(A,B;Z)$.
\begin{definition}
An operation\footnote{For our purposes in this paper the definition
is a scaled back version of the one in \cite{HWZ3}.} for $\pi$
consists of the following data:
\begin{itemize}
\item[1)] A degeneration structure $(S,R)$, where $S=\pi_0(X)$
together with an indexing of $X$.
\item[2)] A map $\circ_{(A,B;C)}:[A]\times [B]\rightarrow [A,B;C]$,
increasing the degeneracy by $1$, which is a sc-diffeomorphism so
that the following properties hold:
\begin{itemize}
\item[2.1)] With $(A,B,C)\leftarrow (A,E)\leftarrow (D)$ and
$(A,B,C)\leftarrow (F,C)\leftarrow (D)$ dual degeneration sequences
we have for $a\in[A], b\in[B]$ and $c\in [C]$ that
$$
a\circ_D\circ (b\circ_E c) = (a\circ_F b)\circ_D c.
$$
\item[2.2)] If $x=a\circ_E b = a'\circ_E b'$ with $a\in [A]$,
$a'\in [A']$ and $d(a)=d(a')$ then $A=A'$.
\end{itemize}
\end{itemize}
We also assume that $\circ$ extends to linear isomorphisms on the
corresponding fibers.
\end{definition}
The definition just given is too special for most applications we
have in mind\footnote{For example for SFT we need that the operation
is compatible with morphisms, has to be multi-valued, and rather
than being as in item 2) sc-diffeomorphisms, we would have only
covering maps.} Nevertheless it is the most instructive one. We
outline necessary modifications later on.

 As a consequence of the
axioms of an operation every element $x\in X$ with $d(x)\geq 1$ is
decomposable. Every element has a prime decomposition in the sense
that for $x\in [Z]$ with $d(x)\geq 1$ there is a uniquely determined
sequence $(A_0,...,A_{d(x)})$ so that there are uniquely determined
$a_i\in [A_i]$ and so that following any degeneration sequence
$$
(A_0,...,A_{d(x)})\leftarrow...\leftarrow (Z)
$$
the associated $\circ$-maps map $x$ onto $(a_0,..,a_{d(x)})$. We
define the spectrum of $x$, denoted by $\sigma(x)$ as the
generalized relator
$$
\sigma(x)=(A_0,...,A_{d(x)};Z).
$$
We call a Fredholm section $f$ of $\pi:Y\rightarrow X$ compatible
with the operation provided for $(A,B;C)$ and $a\in[A]$, $b\in
[B]$ we have
$$
f(a\circ_C b) = f(a)\circ_C f(b).
$$
If $f$ is a section of $\pi:Y\rightarrow X$ define a perhaps
multi-valued section $f\circ f$ on $\partial X$ by
$$
(f\circ f)(z) = \{f(a)\circ_C f(b)\ |\ a\circ_C b=z\}.
$$
 Also define $\partial f$ to be the restriction of $f$ to $\partial
X$. Clearly the compatibility with the operations means that
$$
\partial f =f\circ f.
$$
We call this the "Master Equation". If $f$ is
multi-valued\footnote{Multi-valued sections have to be considered in
the polyfold case due to serious transversality issues.},
compatibility would be defined by the same equation. If $Q$ is any
subset of $X$ we can define
$$
\partial Q= Q\cap(\partial X),
$$
where we recall that $\partial X=\{x\in X\ |\ d(x)\geq 1\}$, and
$$
Q\circ Q=\{a\circ_C b\ |\ (A,B;C)\in R,\ a\in A\cap Q,\ b\in B\cap
Q\}.
$$
Then we can define compatibility of $Q$ with the operation by
$$
\partial Q = Q\circ Q.
$$
The Master Equation turns up everywhere in the theory. We know
already what is means that $f$ is compatible with $\circ$. If
$K=f^{-1}(0)$ and $f$ is compatible with $\circ$ we have $\partial
K=K\circ K$. Not surprisingly the perturbation theory for a
$\circ$-compatible $f$ has to be consistent with $\circ$ by
requiring that the perturbation $s\in\Gamma^+(\pi)$ satisfies
$$
\partial s =s\circ s,
$$
It turns out that there is a very good abstract $\circ$-compatible
perturbation theory. However, one should point out that depending on
the circumstances the transversality-theory can be quite complicated
and elaborate. This is in particular the case if the indexing
involves so-called diagonal relators, i.e. relators of the form
$(A,A;Z)$. In that case we have a situation where we can take two
copies of the same object and can create a new one. This in general
implies an additional symmetry in the problem which causes some
problems in transversality questions, which in most cases can only
be resolved by multi-valued perturbations. The resulting moduli
spaces then will be branched manifolds, see \cite{McDuff}, rather
than manifolds and counting of solutions can usually be done only
over the rational numbers\footnote{This then in general requires
that the Fredholm problems are orientable.}. Similarly, in the
polyfold context, we need to have the perturbations not only
compatible with the operation, but also with the morphisms. In
general, transversality issues can only be resolved by multi-valued
perturbations.

Let us observe that there are many possible generalizations of the
above definition. For example we might have some group action on the
set $S$ which happens for example in SFT where $\pi_2$ acts on the
connected components. The (special) definition we have given for an
operation does not take into account that we might have different
recipes which allow us to associate to a pair of points a bunch of
other points in different components. Moreover, it might occur that
the map $(a,b)\rightarrow a\circ_C b$ is not a diffeomorphism
between $A\times B$ and a face $F$, but only a finite-to-one
covering map. Then, of course, we would like to incorporate
symmetries which is best dealt with the groupoid set-up we discussed
previously. Clearly, it is not necessary to take $S=\pi_0(X)$. In
fact, for a given problem there might be better criteria for
indexing the space than just connectivity components. For example
one would like to bundle several connected components and  denote
this subspace by a letter $A$ and the collection of all these
subspaces will be the set $S$. We refer the reader to \cite{HWZ3}
for a quite comprehensive picture.

Let us informally illustrate some of the ideas in the case of
Morse-theory. In the case of a Morse-function $\Phi:M\rightarrow
{\mathbb R}$ we can take $S$ to consist of all pairs $(a,b)$ of
critical points with $a<b$. The relators consist of all triple
$((a,b),(b,c);(a,c))$. The operation is defined by associating to
$(a,b)$ the component $[(a,b)]$ which is the union of all
$X(a_0,..,a_k)$ with $a_0=a$ and $b=a_k$. If $(a,b)\in S$ then
$R_{(a,b)}$ consist of all admissible symbols $((a,c),(c,b);(a,b))$
and
$$
[(a,c),(c,b);(a,b)]=[(a,c)]\times [(c,b)].
$$
The map $\circ_{((a,c),(c,b);(a,b))}$ associates to elements $x\in
[(a,c)]$ and $y\in [(c,b)]$ the broken trajectory obtained from
$x$ and $y$.  The spectrum of a broken trajectory $x$ connecting
$a_0,a_1,...,a_n$ is
$$
\sigma(x)=((a_0,a_1),(a_1,a_2),..,(a_{n-1},a_n);(a_0,a_n)).
$$
Of course, we could also take a much finer indexing by, in fact,
taking connected components for the indexing, which can be described
by elements in the relative first homotopy $\pi_1(M;a,b)$.

Coming back to the abstract Fredholm theory it is useful to
introduce an auxiliary concept. In the following, given an indexing
by a degeneration structure $(S,R)$, the symbol $[A]$ denotes the
closed and open subspace of $X$ associated to $A\in S$. Observe that
we allow ourselves to be a little bit more general and do not
require $A$ to be a connected component. Assume we are given a
Fredholm operator with operations $\mathfrak{f}=(\pi,f,\circ)$. For
technical reasons we will restrict to M-polyfolds build on
sc-Hilbert spaces.
 An auxiliary norm $N$ is a continuous map defined on $Y_{0,1}$
with image in $[0,\infty)$ introducing on every fiber $Y_{0,1;x}$ a
complete norm and having some additional properties. First of all
$$
N(h\circ_A k) = \hbox{max}\{N(h),N(k)\}.
$$
Further it is not difficult to see that there is a well-defined
concept of "mixed convergence" for sequences in $Y_{0,1}$, which in
local coordinates would mean that the data in the base is converging
on the $0$-level, whereas the data in the fiber is converging weakly
on level $1$ and strongly on level $0$. We write
$y_k\xrightarrow{\rightharpoonup} y$ for mixed convergence. We
require that $N$ has the property
$$
N(y)\leq \liminf_{k\rightarrow\infty} N(y_k).
$$
One can show that such auxiliary norms exist. There are different
useful concepts of properness for a Fredholm section of
$\pi:Y\rightarrow X$.
\begin{itemize}
\item[1)] We say $f$ is component-proper if on every connected component the
induced operator is proper.

\item[2)] If $f$ is compatible with an operation $\circ$ we say
that it is $\circ$-proper if for every $A\in S$ the operator
induced on $[A]$ is proper.
\end{itemize}
Of course 2) implies 1). Recall that a  subset $K$ of $X$ is said to
be $\circ$-compatible if $\partial K=K\circ K$. Here $\partial
K=(\partial X)\cap K$ and $K\circ K $ is defined by
$$
K\circ K=\{a\circ_C b\ |\ a\in K\cap[A], b\in K\cap[K], (A,B;C)\in
R\}.
$$
For example if $f$ is compatible with an operation $\circ$ then
$f^{-1}(0)$ is  $\circ$-compatible. If $f$ is component-wise proper
so is $f^{-1}(0)$. An open $\circ$-neighborhood $U$ of a
$\circ$-invariant subset $K$ of $X$ is an open set containing $K$
which is $\circ$-invariant. A $\circ^+$-section is a sc$^+$-section
compatible with $\circ$. We denote the whole collection by
$\Gamma^+(\pi,\circ)$,

An important compactness result is the following\footnote{One needs
for the following results that the Fredholm section has the
contraction germ property around every point on level $0$ and not
only around smooth points, see \cite{HWZ3}. In applications the
proof of this type of Fredholm property by exhibiting the
contraction germ form on the $0$-level at non-smooth points is
usually identical to the proofs at smooth points on an arbitrary
level.}
\begin{theorem}
Assume that $\mathfrak{f}$ is a Fredholm operator with operations
which is component-wise proper. Let an auxiliary norm $N$ be given.
Then there exists an open $\circ$-neighborhood $U$ of $f^{-1}(0)$ so
that for every $s\in\Gamma^+(\pi,\circ)$ with support in $U$ and
satisfying $N(s(x))\leq 1$ for all $x$ the section $f+s$ is
$\circ$-compatible and component-wise proper.
\end{theorem}

The next theorem gives an abstract transversality result
\begin{theorem}
Let $\mathfrak{f}$ be a Fredholm operator with operation and $N$ an
auxiliary norm. Assume that $\overline{\mathcal O}$ is the set of
$\circ^+$-sections satisfying $N(s(x))\leq 1$ for all $x$ and having
support in the open $\circ$-neighborhood $U$ of $f^{-1}(0)$. Then
the following holds. There exists a section $s\in {\mathcal O}$ so
that $f+s$ is transversal with the zero-section and in general
position with respect to the boundary strata.
\end{theorem}

{ At this point we would like to count solutions on components where
the Fredholm index is $0$. We can do this over ${\mathbb Z}_2$, or
if we can orient the determinants of the linearized Fredholm
sections in a coherent way, we can work with more general
coefficients. As we pointed out earlier there is a natural class of
strong sc-connections. Taking the covariant derivative of a section
with respect to two such connections the difference will be a
sc$^+$-operator, in particular compact. As a consequence there is a
convex set of possible linearisations differing by compact
operators. We will not address these issues here, but point out that
for SFT, orientations over ${\mathbb Z}$ have to be taken, since due
to intrinsic "orbifold difficulties" counting has to be done over
the rational numbers.}

In our example for the Morse-function $\Phi$ let us do the counting
over ${\mathbb Z}_2$. Given $(a,b)$ with Morse-index difference $1$,
i.e. $m(b)-m(a)=1$, denote by $Q(a,b)$ the number of solutions of
$(f+s)(x)=0$ in $X(a,b)$, otherwise define the map to be $0$.
Observe that the Fredholm index $i(a,b)$  is given by
$$
i(a,b)=m(b)-m(a)-1.
$$
Counting gives a map
$$
Q:S\rightarrow {\mathbb Z}_2
$$
For the convolution product $\ast$ one can verify that $Q\ast Q=0$.
We can partition $S$ into even and odd elements by saying that
$(a,b)$ is odd if the difference of the Morse-indices is odd and
even otherwise\footnote{On the level of Fredholm indices a component
is even (odd) if the Fredholm index is odd (even) since we divided
out by the ${\mathbb R}$-action.}. If we consider the space of all
functions from $S$ to ${\mathbb Z}_2$, say $C(S,{\mathbb Z}_2)$, we
therefore obtain a decomposition
$$
C(S,{\mathbb Z}_2) = C_0(S,{\mathbb Z}_2)\oplus C_1(S,{\mathbb
Z}_2).
$$
Using the convolution product we can define a commutator
compatible with the grading (For two odd elements the commutator
has a $(+)$-sign, for all the other cases it has a $(-)$-sign).
Then
$$
[Q,Q]=2Q\ast Q=0.
$$
At this point we have produced the data $(S,Q)$. The next step is
the representation of this data. There are in fact different
possibilities. For example we can define
$$
D_Q:C(S,{\mathbb Z}_2)\rightarrow C(S,{\mathbb Z}_2)
$$
by
$$
D_Q(\lambda) = [Q,\lambda].
$$
Then $D_Q^2=0$ and we obtain a Homology group. If we do it for $S^2$
with the height function we see that $S$ consists of one point, say
$\ast$, which is even. The counting function $Q$ is $0$. If we take
a different Morse-function on $S^2$ having for example four critical
points we obtain a more complicated homology. This kind of homology
can always be defined for any abstract situation and has some
invariance properties with respect to small perturbations. It seems
that  this homology in our Morse-theory example can be used to
define an isotopy invariant for Morse-functions with prescribed type
of critical points. If we know more about the data we can do some
"representation theorem". For example in our Morse-theory case
denote by $V$ the vector space of maps $Cr(\Phi)\rightarrow {\mathbb
Z}_2$. We can use $Q$ to define a linear operator
$$
Q:V\rightarrow V,
$$
by
$$
(Qh)(a) =\sum_{(a,b)\in S} Q(a,b)h(b).
$$
Then $Q^2=0$ and the  homology of $(V,Q)$ is the usual homology with
${\mathbb Z}_2$-coefficients, which would be invariant under
arbitrary (generic) perturbations of $\Phi$ if $M$ is compact.

The following fact should be pointed out. The data $(S,Q)$ is a
homologically invariant way of counting solutions for a Fredholm
problem with operation. This kind of counting or versions applying
more sophisticated topological methods to the whole moduli space
always can be done in the abstract framework. What we have
illustrated above is the instance of a representation theory for
this "counting data". Depending on additional structure of the
counting data we might be able to use it to construct now algebraic
objects. For example in SFT super-Weyl-differential algebras, see
\cite{EGH}. Note, however, that when \cite{EGH} was written, the
technical tools for separating out the analytical, topological and
algebraic aspects as cleanly as it is possible now, had not been
developed.

\subsection{Comments}
There is an interesting algebra of formal differential equations in
the background which can be used to describe homotopies and other
relevant concepts. We will not explain this here and refer the
reader to \cite{HWZ3}. The reader might also have a look at
\cite{EGH} where a certain number of mysterious differential
equations arise in the study of cobordisms and homotopies. These
formulas have a completely rigorous definition within some new
algebra of formal differential equations.

Counting solutions in general calls for a theory of orientations. An
orientation for Fredholm problems is by definition an orientation
for the determinant bundle associated to the point-wise linearized
Fredholm operator. For a general Fredholm problem there is no
orientation. However,  if the linearized operators are compactly
homotopic to complex linear Fredholm operators such a determinant
bundle is orientable. This underlying reason is for example valid in
the case of Gromov-Witten Theory. In SFT the orientation question is
very subtle and can be viewed as an extension of the study of
coherent orientations in Floer Theory, see \cite{FH}. A basic
outline for the general orientation questions in SFT has been given
in \cite{EGH}, and worked out in some variation in
\cite{BourgeoisMohnke}. In our general situation where we deal with
Fredholm Theory with Operations a coherent orientation will refer to
a choice of orientation of the linearized Fredholm section
(Linearized using a special class of connections which is intrinsic
to the theory of strong bundles) so that on every component the
orientation changes continuously, i.e. by prolongation along paths,
and, most importantly for a given relator $(A,B;Z)$ there is a
well-defined relationship between the orientations ${\mathfrak{
o}}_A$, ${\mathfrak{ o}}_B$ and ${\mathfrak{ o}}_Z$. Clearly using
the operation $\circ$ one can come up with a standardized procedure
how to define an orientation for $Z$ if we have given ones for $A$
and $B$. This results in an orientation ${\mathfrak{
o}}_A\circ_Z{\mathfrak{ o}}_A$ which is obtained from orientations
for $A$ and $B$ following some specific constructions (i.e.
conventions). Assuming that $Z$ is already oriented the requirement
for a coherent orientation is the relationship
$$
{\mathfrak{ o}}_A\circ_Z{\mathfrak{ o}}_B
={(-1)}^{p(A)}{\mathfrak{o}}_Z.
$$
Here $p(A)$ is some parity associated to $A$ and related to the
Fredholm index on $A$. As there are important sign conventions in
Homology theory we also have crucial sign conventions here. The
above formula contains such a convention which is (keeping in mind
that counting is a homological process) compatible with the sign
conventions of Homology.

Let us also mention the following. Rather than counting the
solutions in the components where the Fredholm index is $0$ we could
look at the full solution space $K=f^{-1}(0)$ which of course
satisfies
$$
\partial K = K\circ K.
$$
In general, i.e. if study $\circ$-proper Fredholm sections of a
strong polyfold bundle, we need multi-valued $\circ$-compatible
sc$^+$-perturbations. In that case we will obtain a branched
finite-dimensional manifold $M$, see \cite{McDuff}. Many
differential geometric concepts make sense in this context. Of
course we have as a consequence of the $\circ$-compatibility of the
perturbation that $\partial M=M\circ M$. It would be useful to
develop algebraic topology concepts for spaces with an operation. It
is clear that one can develop a theory of differential forms
satisfying $\partial \omega=\omega\circ \omega$, but other concepts
should also carry over.

In the interesting paper \cite{BC} by Barraud and Cornea it is shown
how the equation $\partial A = A\circ A$ can be exploited
algebraically in Floer theory without bubbling, and even in the case
of Morse theory leading to results going far beyond the discussion
of the usual Morse complex. A crucial input in their discussion is
some kind of representation of the moduli spaces in loop spaces. The
moduli spaces have boundaries and they introduce a spectral sequence
which in effect allows them to systematically forget the boundaries
allowing them to define a reduced homology class which is a new
invariant representing  higher dimensional moduli spaces. Their
results depend on a representation theory of the moduli spaces in
loop spaces which is quite natural in their context. It would be
interesting to study the question of representations of other moduli
spaces. Floer theory with bubbling seems already to be  a
challenging start.

Finally let describe in some detail what a general theory of
operations should be. Of course, one would like a notion as general
as possible, subject to a constraint. Namely one wants easy axioms
and within such a general theory one would like to have an abstract
perturbation and transversality theory. Currently we have a theory
which goes beyond what we described here, covers FT, CH and SFT, but
still doesn't achieve what we describe now (but is close). Assume we
are given a polyfold $Z$ with boundary $\partial Z$. One would like
to explain $\partial Z$ in terms of  (fibered) products of its
components. So given a pair of points $(a,b)$ for suitable
components $A$ and $B$ one can produce a points $z$ in a suitable
face of $\partial Z$, i.e. one has a relator $(a,b;z)$. If we now
vary $a$ and $b$ the target $z$ should change smoothly in dependence
of $a$ and $b$. Let us view this, and refer to it, as a "smoothly
changing recipe". Now we could imagine that we have a whole family
of smoothly changing recipes. Then for two points $(a,b)$ there
might indeed be several relators $(a,b;z_j)$. On the other hand
sometimes different recipes applied to different points might imply
the same result. With other words one would like to have a theory of
families of interacting smoothly changing recipes described by a
simple set of axioms, so that the concept allows to develop an
abstract transversality and perturbation theory. Of course,
everything should be so general that the concrete theories of
interest fit into this scheme, but in addition, the level of
generality is right in the sense that its description and necessary
constructions are relatively easy.

\section{Gromov-Witten Theory}
We begin with Gromov-Witten theory. The scale-analysis which has to
be carried out to construct the ambient spaces of SFT is not much
more difficult than that needed for Gromov-Witten theory. Besides
that it is interesting to have a polyfold set-up in that case as
well.

An important input in constructing the polyfold set-up for GW or SFT
is the Deligne-Mumford theory of stable Riemann surfaces, however in
a modified form. Our description of Deligne-Mumford theory and its
modifications geared towards applications in SFT is taken from
\cite{EHWZ1}.

\subsection{Deligne Mumford Type Spaces} For the
analysis of SFT it is important to understand certain variants of
the Deligne-Mumford theory of stable Riemann surfaces. In fact there
are a certain number of issues which will be important to understand
and which are not classical and nonstandard and deal with the fact
that the (smooth) SFT-constructions need differentiable structures
which are (with exceptions) not compatible with the standard
DM-Theory. The DM-background as used in SFT is being developed in
much detail in \cite{EHWZ1}. We give some minimal background here.

We consider tuples $(S,j,M,D)$, where $(S,j)$ is a closed Riemann
surface, $M\subset S$ a finite subset of un-numbered points, and $D$
a finite collection of un-ordered pairs of points $\{x,y\}$, where
$x\neq y$. We call $(S,j,M,D)$ a noded Riemann surface with
(un-ordered) marked points. We say it is connected, provided the
topological space obtained by identifying $x$ with $y$ for every
nodal pair, is a connected topological space. Moreover, we assume
that $\{x,y\}\cap \{x',y'\}\neq \emptyset$ implies
$\{x,y\}=\{x',y'\}$. In addition the union $|D|$ of all these
two-point sets is disjoint from the points in $M$. We will refer to
$M$ as the marked points and $D$ the set of nodal pairs. We say that
$(S,j,M,D)$ is equivalent to $(S',j',M',D')$ provided there exists a
biholomorphic map $\phi:(S,j)\rightarrow (S',j')$ such that
$\phi(M)=M'$ and $\phi(D)=D'$, where
$$
\phi(D):=\left\{\left\{\phi(x),\phi(y)\right\}\
|\left\{x,y\right\}\in D\right\}.
$$
We call $(S,j,M,D)$ stable provided its automorphism group is
finite. Denote by $\overline{\mathcal N}$ the collection of all
equivalence classes of stable noded Riemann surfaces with marked
points. In SFT one will need usually somewhat more complicated
objects. For example one needs some of the points in $M$ to be
ordered, some un-ordered, and some carrying a distinguished oriented
real line in their tangent space. However stripping away all
additional data we end up with the objects just introduced. Given an
equivalence class $\alpha=[S,j,M,D]\in \overline{\mathcal N}$ we can
define its type as follows. We associate to $\alpha$ the isomorphism
class of a decorated graph by declaring the components of $S$ to be
the vertices together with a number giving the genus and a second
number giving the number of points from $M$ lying on the component.
In addition we draw an edge for every $\{x,y\}$ connecting the
components on which the points $x,y$ reside, see Figure
\ref{Figure4}. Note that $x,y$ can lie on the same component. Let us
define the arithmetic genus $g_a$ of $\alpha$ by
$$
g_a = 1+ \sharp D +\sum_{C} [g(C)-1]
$$
where the sum is taken over all connected components of $S$. Assume
that $\alpha$ is stable, i.e. for every vertex the union of twice
$g(C)$ plus the number of marked points plus the number of nodal
points should be at least $3$.

Denote by $\Omega^1(\alpha)$ the space of smooth maps associating to
$x\in S$ a complex anti-linear map $\phi(x):(T_xS,j)\rightarrow
(T_xS,j)$. Let $\Gamma(\alpha)$ be the space of vector fields on $S$
which vanish at the points in $D$ and $M$. Then the Cauchy-Riemann
operator $\bar{\partial}$ defines a Fredholm map
$$
\Gamma(\alpha)\rightarrow \Omega^{0,1}
$$
which is injective and has index $3-3g_a-\sharp M-\sharp D$. We
denote by $H^1(\alpha)$ the complex vector space of dimension
$3g_a+\sharp M+\sharp D -3$ defined by
$$
H^1(\alpha) = \Omega^{0,1}/\bar{\partial}\Gamma(\alpha).
$$
If we denote by $\tau$ the type of $\alpha$ and by
$\overline{\mathcal N}_{\tau}$ the subset of $\overline{\mathcal N}$
of type $\tau$ then one can identify $H^1(\alpha)$ with the (orbi-)
tangent space of $\overline{\mathcal N}_{\tau}$ at $[\alpha]$. Note
that the automorphism group $G$ of $\alpha$ acts on $H^1(\alpha)$
and that isomorphisms $\phi:\alpha\rightarrow\alpha'$ induce
isomorphisms $H^1(\alpha)\rightarrow H^1(\alpha')$. Given a smooth
family of complex structures $v\rightarrow j(v)$ on $S$, say with
$v\in V\subset E$, where $E$ is some finite-dimensional vector space
with $j(0)=j$ we can take its derivative $Dj(v)$ at $v$ and observe
that it induces a linear map
$$
E\rightarrow H^1(\alpha_v)
$$
where $\alpha_v=(S,j(v),M,D)$, called the Kodaira differential and
denoted by $[Dj(v)]$.

It is convenient to take $E=H^1(\alpha)$ which is a complex vector
space with a natural action of the automorphism group $G$ of
$\alpha$ on it. Let $V\subset H^1(\alpha)$ be a $G$-invariant open
neighborhood of $0$. We call a family $v\rightarrow j(v)$ with
$j(0)=j$ effective if at every $v\in V$ the Kodaira differential is
a real linear isomorphism. We call it complex if at every $v$ the
Kodaira differential is complex linear. We call it symmetric
provided for every $v \in V$ and $g\in G$ the map
$$
g:(S,j(v),M,D)\rightarrow (S,j(g\ast v),M,D)
$$
is an isomorphism. A family $v\rightarrow j(v)$ is called good if it
is effective, symmetric and smooth. One can show that for every
$\alpha$ there exists a good complex family so that the maps
$$
v\rightarrow [S,j(v),M,D]
$$
define a smooth orbifold structure for $\overline{\mathcal
N}_{\tau}$ with holomorphic transition maps. Next we would like to
address the problem of defining smooth orbifold structure on
$\overline{\mathcal N}$. This is more delicate if one keeps in mind
an important goal: Such smooth structures should be compatible with
a "to be constructed" theory of smooth structures for the moduli
spaces (or their perturbations) of stable finite energy surfaces, a
notion we will introduce later. For the moment it suffices to know
that these are the solutions of our Fredholm problems arising in
SFT.

Start with $\alpha$ of type $\tau$. Then one can find a good complex
family $v\rightarrow j(v)$ with $j(v)=j$ near the points in $D$ and
$M$. A small disk structure consists of a family of disks $D_x$, for
every nodal point $x$, with smooth boundaries so that their union is
invariant under the $G$-action. Assume that $j(v)=j$ on these disks.
Fix for every $x$ a biholomorphic map $\bar{h}_x:(D,0)\rightarrow
(D_x,x)$ and complex anti-linear maps
$\varphi_{(x,y)}:T_yS\rightarrow T_xS$ so that
$\varphi_{(x,y)}^{-1}=\varphi_{(y,x)}$ and the following
compatibility holds
$$
T\bar{h}_y(0)^{-1}\circ\varphi_{(y,x)}\circ T\bar{h}_x(0):{\mathbb
C}\rightarrow {\mathbb C}
$$
is complex conjugation. Next take for every nodal pair $\{x,y\}$ a
copy of the complex plane ${\mathbb C}_{(x,y)}$ and let
$$
N=\bigoplus_{\{x,y\}\in D} {\mathbb C}_{(x,y)}.
$$
Then there is a natural unitary action of $G$ on $N$ which will be
compatible with the following construction where small elements in
$N$ occur as gluing parameters. Pick $\{x,y\}$ and define
$h_x(s,t)=\bar{h}_x(e^{-2\pi(s+it)})$ and
$h_y(s',t')=\bar{h}_y(e^{2\pi(s'+it')})$, where $t,t'\in S^1$ and
$s\leq 0$ and $s'\geq 0$. Note here that we pick for one of the
points positive holomorphic polar coordinates and for the other
negative ones. However the following construction does not depend on
this choice. For $R\geq 0$ define $D_x^R$ and $D_y^R$ by
\begin{eqnarray*}
& D_x^R=\{z\in D_x\ |\ z=h_x(s,t),\ s\in [0,R],\ t\in S^1\}&\\
& D_y^R=\{z\in D_y\ |\ z=h_y(s',t')\ s'\in [-R,0],\ t\in S^1\}&.
\end{eqnarray*}
Here we arrive at the important point where we have to discuss the
gluing (also in this context sometimes called plumbing).  Observe
that we have to make a choice how a gluing parameter
$a=a_{\{x,y\}}\in {\mathbb C}_{\{x,y\}}$ is being converted into a
gluing length $R$ and a gluing angle $\vartheta$. There is not too
much choice with the angle, but there are plenty of choices for $R$.
Observe that this discussion parallels the in the Morse-theory
example. Let us recall the notion of a gluing profile:
\begin{definition}
A gluing profile $\varphi$ is a diffeomorphism
$$
\varphi:(0,1]\rightarrow [0,\infty).
$$
\end{definition}
Fix a gluing profile $\varphi$ and define if $0<|a_{\{x,y\}}|< 1$
$$
R=\varphi(|a_{\{x,y\}}|)\ \hbox{and} \
a_{\{x,y\}}=|a_{\{x,y\}}|e^{-2\pi i\vartheta}.
$$
Now call points $z\in D_x^R$ and $z'\in D_y^R$ equivalent provided
$$
s-s'=R\ \hbox{and}\  t-t'=\vartheta \ (\hbox{mod}\ 1).
$$
Doing this for all gluing parameters $a$ we obtain a new family of
Riemann surfaces $\alpha_{(v,a)}=(S_a,j(v,a),M_a,D_a)$. Moreover
every $g\in G$ defines in a natural way an isomorphism
$$
g:\alpha_{(v,a)}\rightarrow \alpha_{g\ast(v,a)}.
$$
We have
\begin{theorem}
The space $\overline{\mathcal N}$ possesses a natural paracompact
topology which is locally compact. The following holds:
\begin{itemize}
\item Fixing a gluing profile $\varphi$ and taking for every
$\alpha$ of type $\tau$ a good  family  the above construction
defines via
$$
(v,a)\rightarrow [\alpha_{(v,a)}],
$$
if restricted to a small enough $G$-invariant neighborhood $(0,0)$ a
family of $C^0$-uniformizers for $\overline{\mathcal N}$. \item
Taking the gluing profile $\varphi(x)=-\frac{1}{2\pi}\ln(|x|)$ and
starting with good complex families the associated uniformizer
define a holomorphic orbifold structure.

\item Starting with the gluing profile $\varphi(x)
=e^{\frac{1}{x}}-e$ and any good family the construction gives a
smooth orbifold structure together with a natural homotopy class of
almost complex structures.
\end{itemize}
\end{theorem}
This result is proved in \cite{EHWZ1}. The logarithmic gluing
profile gives in fact the classical Deligne-Mumford structure.
Unfortunately it is not compatible with a smoothness theory for SFT
or even with Gromov-Witten theory in the smooth non-integrable case.
The exponential gluing profile however works very well. The identity
map from the exponential smooth structure to the DM-structure is
however smooth. The above theorem, and in particular some estimates
which one obtains during the proofs are important for the
SFT-theory.

Let us give a somewhat different description of the previous
results, which where formulated in the classical V-manifold or
orbifold language, by describing them via the groupoid approach (see
\cite{MoerMr} for the relevant Lie groupoid theory). The space
$\overline{\mathcal N}$ constructed above will be viewed as the
orbit space of a Lie groupoid. To motivate the "groupoid approach"
consider the huge category (not a set) where the objects are tuples
$\alpha:=(S,j,M,D)$ of stable Riemann surfaces. Let us define a
morphism $\alpha\rightarrow\alpha'$ to be a biholomorphic map
$$
\phi:(S,j)\rightarrow (S',j')
$$
satisfying $\phi(D)=D$ and $\phi(M)=M'$. Clearly, every morphism is
invertible. Calling two objects equivalent if there is a morphism
between them, we can build the quotient of the category, which is a
set and, in fact, the orbit space we are interested in. Since the
quotient space is a set there is an enormous redundancy in the
original category and we have to cut it down to a suitable full
subcategory which still has the same quotient space, but is a set.
In order to do so in a sensible way observe that it is possible to
say in our category that two objects which are not isomorphic are
close (up to isomorphism), see the discussion of the DM-theory in
\cite{BEHWZ}. This then allows to put a topology on the orbit space.
Having this all in place we notice that we can describe all
equivalence classes by a suitable set of parameterized models. Pick
an element $\alpha=(S,j,M,D)$ and fix a good family $v\rightarrow
j(v)$, together with a small disk structure and define
$$
(a,v)\rightarrow \alpha_{(a,v)}=(S_a,j(a,v),M_a,D_a).
$$
If we restrict this map to a sufficiently small open neighborhood of
$(0,0)\in N\times E$, say $W$, it would induce a unifomizer as
constructed before. Consider the set
$$
\Delta:=\{(a,v,\alpha_{(a,v)})\ |\ (a,v)\in W\}.
$$
Then $\Delta$ possesses a natural smooth manifold structure by
requiring that the map
$$
\Delta\rightarrow W:((a,v),\alpha_{(a,v)})\rightarrow (a,v)
$$
is smooth. Construct in a similar way $\Delta'$. Assume that we have
an isomorphism
$$
\phi:\alpha_{(a,v)}\rightarrow \alpha'_{(a',v')}.
$$
Then one can show, that $\phi$ lies in a uniquely determined family
of isomorphisms $\phi_{(b,w)}$ between $\alpha_{(b,w)}$ and
$\alpha'_{u(b,w)}$.  Here $u$ stands for a uniquely determined
smooth local diffeomorphism
$$
(b,w)\rightarrow (b',w')=u(b,w)
$$
with $u(a,v)=(a',v')$. We view
$$
(u(b,w),\phi_{(b,w)},(b,w)):((b,w),\alpha_{(b,w)})\rightarrow
(u(b,w),\alpha'_{u(b,w)})
$$
as a morphism. Note that the set of morphisms
$$
\Delta\rightarrow \Delta'
$$
has the structure of a smooth manifold of the same dimension as
$\Delta$. We can now find a countable family of $\Delta$'s, say
$\Delta_i$, $i\in I$, so that with $X$ defined by the disjoint union
$$
X=\coprod_{i\in I}\Delta_i
$$
the map $X\rightarrow \overline{\mathcal N}$ which assigns to an
element $(a,v,\alpha_{(a,v)})$ the isomorphism class of
$\alpha_{(a,v)}$ is a surjection. The previous discussion already
identified the morphisms and we know that the collection of all
morphisms ${\bf X}$ has a manifold structure as well. It follows
immediately from the construction that the associated category
$\mathfrak{X}$ is an \'etale Lie groupoid\footnote{We could modify
it and make it even proper.} whose orbit space is
$\overline{\mathcal N}$. Consider $(\mathfrak{X},\beta)$, where
$$
\beta:|\mathfrak{X}|\rightarrow \overline{\mathcal N}
$$
is the obvious homeomorphism. Then  we can say that the pair defines
an orbifold structure on $\overline{\mathcal N}$. As in the polyfold
case we can define equivalence of two such pairs by refinement. This
gives an alternative way of defining an orbifold structure. The
method will immediately generalize to the case where we also have a
map on $S$. This applies then to Gromov-Witten theory.

\subsection{The Splicings for Gromov-Witten Theory}
The main ingredient, as in the Morse-Theory case, is a gluing and
anti-gluing procedure. The situation in GW differs somewhat from the
latter case, since the isolated critical points in Morse-Theory are
replaced by constant loops, i.e. the symplectic manifold. With other
words we are dealing in fact with a Morse-Bott situation. This
requires some modifications.

We begin with the discussion of the procedure how maps $u$ on a
nodal Riemann surface $S$ can be used to describe neighboring curves
on (un-noded or glued) Riemann surfaces $S_a$ by the so-called
gluing construction $\oplus_a(u):S_a\rightarrow W$. This
construction which can be carried out in local coordinates and
implanted into a manifold setting by a chart will be referred to as
"nonlinear gluing". The gluing also exists on the level of vector
fields and leads to vector fields along glued maps. As in the Morse
theory case there is a certain ambiguity in gluing in the sense that
gluing of different maps might lead to the same result. Here is the
place where anti-gluing comes in. It is defined on the level of
vector fields and precisely resolves this ambiguity. Moreover, the
combination of gluing and anti-gluing, called total gluing, leads to
the abstract concept of splicing which is crucial for our theory.
With other words the whole philosophy is as discussed in the Morse
theory case. Note, however, that there will be some modifications.
This is mainly due to the fact that we deal with a Morse-Bott,
rather than a Morse-situation.

We begin with the basic gluing $\oplus$. Fix a smooth cut-off
function
$$
\beta:{\mathbb R}\rightarrow [0,1]
$$
so that
\begin{eqnarray*}
&\beta(s)+\beta(-s)=1\ \hbox{for all}\ \ s\in {\mathbb R}&\\
&\beta'(s)<0\ \ \hbox{for all}\ \ s\in (-1,1).&
\end{eqnarray*}
Also recall our gluing profile $\varphi:(0,1]\rightarrow [0,\infty)$
defined by
$$
\varphi(r) = e^{\frac{1}{r}}-e.
$$
Let us first assume that we are given two abstract disk-like Riemann
surfaces $D_x$ and $D_y$ with smooth boundaries and interior points
$x$ and $y$. We view $\{x,y\}$ as a nodal pair. Hence $(D_x\cup
D_y,\{x,y\})$ is a noded Riemann surface. As explained in the
background chapter about DM-theory we can glue the two disk given a
gluing parameter. We will however keep a little bit more
information. Let us also assume that we have given holomorphic polar
coordinates centered around $x$ and $y$. To be more precise we are
given biholomorphic maps $\bar{h}_x:(D_x,x)\rightarrow (D,0)$ and
$\bar{h}_y:(D_y,y)\rightarrow (D,0)$. Then define
$$
h_x:{\mathbb R}^+\times S^1\rightarrow D_x:(s,t)\rightarrow
\bar{h}_x(e^{-(2\pi(s+it))})
$$
and
$$
h_y:{\mathbb R}^-\times S^1\rightarrow D_y:(s',t')\rightarrow
\bar{h}_y(e^{2\pi(s+it)}).
$$
If $a\in {\mathbb C}$ is a complex number with $|a|\leq \frac{1}{2}$
we define $Z_0=(D_x\cup D_y,\{x,y\})$, which is our original noded
disk, and if $0<|a|\leq \frac{1}{2}$ we define the cylinder $Z_a$ by
identifying the points $z=h_x(s,t)\in D_x$ for $(s,t)\in [0,R]\times
S^1$ with $z'=h_y(s',t')\in D_y$, where $s-s'=R$ and
$t-t'=\vartheta$ and $(R,\vartheta)$ is associated to $a$. This
cylinder has two preferred sets of coordinates
$$
[0,R]\times S^1\rightarrow Z_a:(s,t)\rightarrow [h_x(s,t)]
$$
and
$$
[-R,0]\times S^1\rightarrow Z_a:(s',t')\rightarrow [h_y(s',t')].
$$
Here $[.]$ means the passing to equivalence classes. We will
abbreviate the first by $[s,t]$ and the second by $[s',t']'.$ Then
$$
[s,t]=[s-R,t-\vartheta]'.
$$
Here comes an important observation. We can define another Riemann
surface as follows. Namely glue $D_x$ and $D_y$ with the same
identification to obtain $\Sigma_a$ which is a simply connected
Riemann surface having two distinguished points $x$ and $y$ and a
distinct annular subregion, namely $Z_a$. This, of course only holds
for $a\neq 0$. If $a=0$ we put $\Sigma_0=\emptyset$.

The purpose of gluing is to associate to a map $u$ defined on
$D_x\cup D_y$, with matching condition at the nodal pair $\{x,y\}$,
a map $\oplus_a(u)$ on $Z_a$. We assume that the target manifold is
${\mathbb R}^{2n}$. Consider continuous maps $u:D_x\cup
D_y\rightarrow {\mathbb R}^{2n}$ with matching condition
$u(x)=u(y)$. We will define a map $\oplus_a(u):Z_a\rightarrow
{\mathbb R}^{2n}$ which we might view as the "nonlinear gluing" of
two maps from a noded Riemann surface with image in the manifold
${\mathbb R}^{2n}$.

Consider the Banach space consisting of $u$ with matching condition
at $\{x,y\}$ so that $u(x)=u(y)=c$ with
$$
u^\pm-c\in H^{3,\delta_0}({\mathbb R}^{\pm}\times S^1,{\mathbb
R}^{2n}).
$$
Here $(u^+,u^-)$ is defined by
$$
u^+=u\circ h_x(s,t)\ \hbox{and}\ \ u^-=u\circ h_y(s',t').
$$
We call $c$ the (common) asymptotic constant of $u^+$ and $u^-$.
Here $\delta_0\in (0,2\pi)$. We equip $E$ with the sc-structure
where level $m$ corresponds to regularity $(m+3,\delta_m)$ and
$(\delta_m)$ is a strictly increasing sequence strictly bounded by
$2\pi$. We define $\oplus_a(u)$ as follows. We put $\oplus_0(u)=u$.
If $0<|a|\leq \frac{1}{2}$ then
$$
\oplus_a(u):Z_a\rightarrow {\mathbb R}^{2n}
$$
given by
\begin{eqnarray*}
&&\oplus_a(u^+,u^-)([s,t])\\
&=&\beta(s-\frac{R}{2})u^+(s,t)\\
&&+(1-\beta(s-\frac{R}{2}))u^-(s-R,t-\vartheta).
\end{eqnarray*}
The interpretation of the procedure so far should be that of a
"nonlinear gluing".

Assume next that $h$ is a vector field along $u$, i.e. $h(z)\in
T_{u(z)}{\mathbb R}^{2n}={\mathbb R}^{2n}$. Again we have a matching
condition $h(x)=h(y)$. Moreover, we define the gluing $\oplus_a(h)$
by the same formula and consider elements $h$ of the same regularity
as that of $u$. Observe that $u+h$ can be interpreted as
$$
z\rightarrow (\exp_u(h))(z)=\exp_{u(z)}(h(z))
$$
for the standard metric. The interpretation of $\oplus_a(h)$ has to
be that of a vector field along $\oplus_a(u)$. Clearly,
$$
\exp_{\oplus_a(u)}(\oplus_a(h))=\oplus_a(\exp_u(h)).
$$

Next we introduce the anti-gluing which is only defined for vector
fields.
 We define the anti-gluing with respect to a reference map
$u:D_x\cup D_y\rightarrow {\mathbb R}^{2n}$ as follows.  Let $c$ be
the asymptotic constant of $u$, without loss of generality say
$c=0$. Then
$$
\ominus_a(h):\Sigma_a\rightarrow T_0{\mathbb R}^{2n}
$$
is given for $a\neq 0$ by
\begin{eqnarray*}
&&\ominus_a(h)([s,t])\\
&=&-(1-\beta(s-\frac{R}{2}))(h^+(s,t)-av_R(h))\\
&&+\beta(s-\frac{R}{2})(h^-(s-R,t-\vartheta)-av_R(h))
\end{eqnarray*}
where
$$
av_R(h)=\frac{1}{2}\left(\int_{S^1}
(h^+(\frac{R}{2},t)+h^-(-\frac{R}{2},t))dt\right).
$$
Moreover $\ominus_0(h)=0$. Observe that
$$
\ominus_a(h)(x)=-\ominus_a(h)(y).
$$
 If we view the above construction as a
description of curves near a noded curve with $u(x)=u(y)=0$, then
one should interpret the image of $\ominus_a(h)$ as $T_{0}{\mathbb
R}^{2n}={\mathbb R}^{2n}$.

The crucial observation of the total gluing construction is that the
map
$$
h\rightarrow \boxplus_a(h)=(\oplus_a(h),\ominus_a(h))
$$
for fixed $a\neq 0$ is a bijective linear sc-operator
$$
E\rightarrow H^3(Z_a,{\mathbb R}^{2n})\oplus
H^{3,\delta_0}_{c}(\Sigma_a,{\mathbb R}^{2n})
$$
where the target spaces are equipped with the sc-structure where the
first component is of class $m+3$ and the second $(m+3,\delta_m)$.
The subscript $c$ stands for antipodal asymptotic constant. Given a
map $w:Z_a\rightarrow {\mathbb R}^{2n}$ near $\oplus_a(u)$ we find a
unique vector field $\eta$ along $\oplus_a(u)$ with
$$
exp_{\oplus_a(u)}(\eta)=w.
$$
Then we find a unique $h$ with
$$
\oplus_a(h)=\eta\ \hbox{and}\ \ \ominus_a(h)=0.
$$
We can define for $|a|\leq \frac{1}{2}$ a sc-decomposition of $E$ by
$$
E=\ker(\oplus_a)\oplus_{sc}\ker(\ominus_a)
$$
and denote by $\pi_a:E\rightarrow E$ the projection onto
$\ker(\ominus_a)$ along $\ker(\oplus_a)$. We will show
\begin{theorem}
The triple ${\mathcal S}=(\pi,E,B_{\frac{1}{2}})$ is a sc-smooth
splicing.
\end{theorem}
Observe that by construction the map
$$
h\rightarrow \oplus_a(\exp_u(h))
$$
where $h\in K_a^{\mathcal S}$ defines a bijection onto the maps
$w:Z_a\rightarrow {\mathbb R}^{2n}$ of class $H^3$. This will be
important in the construction of polyfold charts.

\subsection{The Polyfold Structure for Gromov-Witten Theory}

Consider a compact\footnote{Compactness is not crucial.} symplectic
manifold $(W,\omega)$. We will consider tuple $(S,j,M,D,u)$ where
$(S,j,M,D)$ is a not necessarily stable noded Riemann surface with
an ordered set of marked points $M$. We impose the following
stability condition. For every connected component $C$ of $S$ at
least one of the following holds:
\begin{itemize}
\item[1)] $2g(C)+\sharp(C\cap (M\cup |D|))\geq 3$ or
\item[2)] $\int_C u^\ast\omega>0$.
\end{itemize}
Let us describe the quality of the function $u$ in more detail. We
say that $u$ is of class $H^{m,\varepsilon}$ for some $m\geq 3$ and
$\varepsilon>0$ if $u$ is of class $H^m_{loc}$ on $S\setminus |D|$
and if for every nodal point $x\in \{x,y\}\in D$ there exists a
smooth chart $\varphi$ around $u(x)$ mapping $u(x)$ to $0$ and
positive holomorphic polar coordinates centered around $x$, say
$\sigma:[0,\infty)\times S^1\rightarrow S$ so that
$$
v(s,t)=\varphi\circ u\circ \sigma(s,t)
$$
belongs to $H^{m,\varepsilon}([R_0,\infty)\times S^1,{\mathbb
R}^{2n})$ for some sufficiently large $R_0$. The definition of being
of class $H^{m,\varepsilon}$ does not depend on the choices
involved. Let us remind the reader that
$H^{m,\varepsilon}([R_0,\infty)\times S^1,{\mathbb R}^{2n})$
consists of all functions so that the distributional partial
derivatives up to order  $m$ weighted by $e^{\varepsilon|s|}$ belong
to $L^2$. We call two such tuples equivalent, say
$$
(S,j,M,D,u)\equiv (S',j',M',D',u')
$$
provided there exists a biholomorphic map
$$
\phi:(S,j,M,D)\rightarrow (S',j',M',D')
$$
so that
$$
u'\circ\phi=u.
$$
We have
\begin{lemma}
If $\alpha$ is stable its automorphism group $G$  is finite. In
particular the set of isomorphisms between two stable elements
$\alpha$ and $\alpha'$ is finite.
\end{lemma}
Now fix any sequence $0<\delta_0<\delta_1<...<2\pi$ of increasing
weights and denote by $X$ the collection of equivalence classes
where $u$ is of class $(3,\delta_0)$. We can define a filtration by
nested subsets of $X$ by declaring $X_m$ to consist of all elements
of class $(m+3,\delta_m)$. Let us begin with the topological side of
things.

\begin{theorem}
The space $X$ carries a natural paracompact second countable
topology so that the un-noded curves are dense. Moreover the spaces
$X_m$ carry natural topologies as well, so that every point in $X_m$
has a closed neighborhood which embeds as a precompact set into
$X_{m-1}$ provided $m\geq 1$.
\end{theorem}

Having the defined the topology one can define a polyfold structure
on $X$. More precisely

\begin{theorem}
The second countable paracompact space $X$ defined above admits a
natural polyfold structure.
\end{theorem}

As we will see the construction is quite similar to the groupoid
construction in the Deligne-Mumford case. Natural here means, that
given the filtration $X_m\subset X$ by having fixed the regularity
requirement that level $m$ corresponds to maps of quality
$(m+3,\delta_m)$, and using the exponential gluing profile there is
a standardized way of defining the polyfold structure.

Consider an element $\alpha=(S,j,M,D,u)$ of class $(3,\delta_0)$. By
the stability condition the map $u$ induces on every unstable
component $C$ a non-constant map $u$ with $\int_C u^\ast \omega>0$.
By the Sobolev embedding theorem the map $u$ is of class $C^1$. We
can add a finite  number of (un-ordered) points $\Xi$ to $S$ in
order that $\alpha^\ast:=(S,j,M\cup\Xi,D)$ is a stable Riemann
surface. We view $M^\ast=M\cup \Xi$ as an un-ordered set of marked
points. Assume we are given  a good family $v\rightarrow j(v)$,
where $V\subset E$ is an open invariant neighborhood of $0$ in some
complex space $E$ with an action of the automorphism group $G^\ast$
of $\alpha^\ast$. We find an open neighborhood $O$ in $N$ so that
$$
(a,v)\rightarrow (S_a,j(a,v),M^\ast_a,D_a)
$$
induces a smooth uniformizer for $(a,v)\in O\times V$. Let us assume
now that the set of points $\Xi$ has the property that $\Xi$ is
invariant under $G$ and that given points $z_1,z_2$ belonging to
different orbits the values $u(z_1)$ and $u(z_2)$ are different and
$Tu(z)$ is injective for any such point in $\Xi$. We also assume
that the image $u(\Xi)$ does not intersect $u(|D|\cup M)$.

Having fixed $\Xi$ with these properties pick $z_1,...,z_k\in\Xi$ so
that their orbits are disjoint and the union of the orbits is $\Xi$.
For every orbit $z_i$ pick a complement $C_i$ of the image of
$Tu(z_i)$. Let us define $w_i=u(z_i)$ for $i=1,...,k$ and let us
denote by $w_{k+1},...,w_{k+\ell}$ an enumeration of the points in
$u(|D|)$. We fix diffeomorphic charts around $w_1,...,w_{k+\ell}$,
by
$$
\varphi_i:({\bf R}(w_i),w_i)\rightarrow ({\mathbb R}^{2n},0).
$$
Here the closures of the domains are mutually disjoint. Let us
denote by ${\bf R}_r(w_i)$ the preimage of the $r$-ball. We pick a
small disk structure for $(S,j,M^\ast,D)$ so that the images of the
$D_x$ are contained in $\bigcup_{i=1}^{k+\ell} {\bf R}_1(w_i)$. Now
we can implant the nonlinear gluing construction so that it is
applicable to maps $u'$ which map the $D_x$ into
$\bigcup_{i=1}^{k+\ell} {\bf R}_4(w_i)$.  Next we consider the space
of $H^{3,\delta_0}$-sections of $u^\ast TW$. By the Sobolev
embedding theorem these sections belong to
$C^1_{loc}(S\setminus|D|)$. We are interested in the subspace of
those sections which over the points in $\Xi$ which belong to $Gz_i$
belong to $C_i$. Note that $C_i$ is a complement for the image of
$Tu(z)$, if $z\in Gz_i$. Let us denote this space of sections by
$Z_{\alpha}$. Consider now the exponential map for a Riemannian
metric which on ${\bf R}(w_i)$ is the pull-back of the standard
metric  by $\varphi_i$.  Then we consider for $(a,v)\in O\times V$
and $\eta$ sufficiently small in $C^0$ and belonging to $Z_\alpha$
the map
$$
(a,v,\eta)\rightarrow
\Phi(a,v,\eta):=(S_a,j(a,v),M_a,D_a,\oplus_a(\exp_u(\eta))).
$$
Recall that for every nodal point we can implant using the charts
the (local) splicing construction giving a splicing
$(\pi,Z_\alpha,U)$, where $U$ is an open neighborhood of $(0,0)\in
O\times V$. Note that $\pi$ only depends on the first component of
$(a,v)$. One could show that the map
$$
(a,v,\eta)\rightarrow [\Phi(a,v,\eta)],
$$
where we pass to equivalence classes, restricted to a sufficiently
small open neighborhood of $0$ in the splicing core, can be viewed
as a uniformizer. Alternatively, as in the Deligne-Mumford case we
define
$$
\Delta=\{(a,v,\eta,\Phi(a,v,\eta)) \ |\ (a,v,\eta)\in P\},
$$
where $P$ is a sufficiently small open neighborhood of $0$ in the
splicing core. Then we can take $\Delta$ and $\Delta'$ and show that
isomorphisms occur in families parameterized by the data in the
domain, i.e. $(a,v,\eta)$. After this, similarly as in the
Deligne-Mumford case we can define a polyfold structure for $X$.
 We
can also define a strong bundle $Y\rightarrow X$. The elements of
$Y$ are equivalence classes of tuples $(S,j,M,D,u,h)$, where
$h:(T_zS,j)\rightarrow (T_{u(z)}W,J)$ is complex anti-linear of
Sobolev regularity $H^{2}_{loc}$ with a particular behavior near the
nodes. The equivalence classes are defined similarly as the
equivalence classes defining $X$. This defines the strong bundle
$$
\pi:Y\rightarrow X.
$$
Then we can define a section $f$ by
$$
f([S,j,M,D,u])=[S,j,M,D,u,\bar{\partial}_{J,j}(u)].
$$
This is, of course, the Cauchy-Riemann operator in our context
extended in the obvious way to nodal surfaces. Recall that
$Y\rightarrow X$ has the structure of a polyfold bundle. In
particular the polyfold structure is represented by a strong
polyfold bundle over some ep-polyfold groupoid. As it turns out, for
suitable representatives the section $f$ is the map induced between
orbit spaces coming from a sc-smooth section functor. Using this we
show in \cite{HWZ2}
\begin{theorem}
The section $f$ of $\pi:Y\rightarrow X$ is a sc-smooth  Fredholm
section.
\end{theorem}
As a consequence abstract perturbation theory using sc$^+$-
multi\-sections is applicable. For such a generic perturbation by a
multi-section the solution space is a locally compact, smooth
branched, weighted manifold, which restricted to every connected
component of $X$, is compact. We can define Gromov-Witten invariants
by integrating suitable quantities over these moduli spaces. One can
study the question of natural homotopy classes  of smooth almost
complex structures for these moduli spaces and many more questions.
In our set-up quite a number of constructions in the GW-theory of
smooth symplectic manifolds become easier than with other
technology.
\section{Symplectic Field Theory}\label{section5}
In this section we explain how SFT fits into the picture of Fredholm
theory with operations.

One of the ideas behind SFT is that it provides tools for computing
Gromov-Witten invariants. It can be viewed as a theory of relative
Gromov-Witten invariants and allows computations by cutting a
compact symplectic manifold along suitable (real) hypersurfaces and
to obtain invariants for the parts taking values in some algebraic
object associated to the hypersurface. It turns out that these
hypersurfaces need some geometrical properties in order to make this
program possible. Then, of course, later one wants to cut along
hypersurfaces in the hypersurfaces and so on. One also would like to
allow  hypersurfaces with singularities, for example chopping up the
symplectic manifold via a triangulation. It seems for example, that
our polyfold language (at least in some suitable
 generalization) is able to deal with the analytic intricacies
arising by degenerating almost complex structures along a
triangulation.

In the above program it turns out that it is important to
understand contact manifolds. Let $(M,\xi)$ be a compact contact
manifold of dimension $2n-1$ (without boundary). We assume that
the contact form is co-oriented. Hence there is a distinguished
class of contact forms inducing $\xi$. Any two of them are related
by multiplication with a positive function. Having fixed such a
$\lambda$ denote its Reeb vector field by $X$. Recall that it is
defined by
$$
\lambda(X)\equiv 1\ \ \hbox{and}\ \ \ d\lambda(X,.)\equiv 0.
$$
{ We like to construct invariants for $(M,\xi)$. We aim at
associating to $(M,\xi)$ an object $O(M,\xi)$ in some category
$\mathfrak{O}$. The construction of $O(M,\xi)$ will be accomplished
in the following way. For $(M,\xi)$ we pick a generic contact form
$\lambda$ inducing $\xi$ and a complex multiplication
 $J:\xi\rightarrow \xi$ compatible with $d\lambda$. Then
 we associate to $(\lambda,J)$ a Fredholm problem with operations
 say $\mathfrak{f}_{(\lambda,J)}$. Usually these problems
 are not generic for any choice of (geometric) data. Then an abstract perturbation
 by sc$^+$-sections (in fact they have to be multi-sections)
will make the data generic and we can produce via counting of
solutions the data for an associated Homology independent of the
(abstract) perturbation. Having made different choices
 $(\lambda,J)$ and $(\lambda',J')$ and taking a generic homotopy
 there is a new associated Fredholm problem $F$ which has a "module structure"
 as left module over $\mathfrak{f}_{(J,\lambda)}$ and as
 right-module over the other. In addition it has some more
 properties which allows to show that taking different
 perturbations leads to the same results if one counts correctly.
We will not discuss this further, but refer the reader to
\cite{HWZ3}.
 We just mention some suggestive formulas which occur everywhere
 in the theory and catch what we mean that the Fredholm problem
 $F$ has a module structure over $f$ and $f'$ as left- and right
 module, respectively:
 $$
 F(x\circ_{B}^L a) = f(x)\circ_{B}^L F(a)\ \hbox{and}
 \ F(b\circ_A^R z) =
 F(b)\circ_A^R f'(z).
 $$
} Of course there is a certain set of axioms "regulating" the
interplay between the "algebra" of the left- and right-
degeneration structures and the way they operate on the middle
problem.
\subsection{Background Material}
We  explain the necessary background material.
\subsubsection{Finite Energy Maps} In order to describe the moduli spaces we are interested in, we
fix a non-degenerate contact form $\lambda$ inducing the contact
structure $\xi$ and respecting the co-orientation. We denote the
associated Reeb vector field by $X$. A contact form is said to be
non-degenerate provided all its periodic orbits, i.e. those of its
Reeb vector field, are non-degenerate.  Recall that this means that
the eigenvalues of the linearized Poincar\'e section maps contain no
root of unity. Non-degeneracy is generic, see
\cite{HWZ-strictly-convex,HWZ2003}. Next we pick a compatible
complex multiplication $J:\xi\rightarrow\xi$ on the contact planes
$\xi$. Compatibility here means that for $h,k\in\xi_m$, $m\in M$,
the map
$$
(h,k)\rightarrow d\lambda(h,J(m)k)
$$
is a positive definite inner product on $\xi_m$. We extend $J$ to
a ${\mathbb R}$-invariant almost complex structure $\tilde{J}$ on
${\mathbb R}\times M$ by defining
$$
\tilde{J}(a,u)(s,h)=(-\lambda(h),J(u)\pi(h) + sX(u)).
$$
Here $\pi:TM\rightarrow \xi$ is the projection along the Reeb
vector field.

We are interested in equivalence classes of solutions of the
nonlinear Cauchy-Riemann associated to $\tilde{J}$. To be more
precise consider tuples\footnote{Of course one also could allow
marked points.}\ $\alpha:=(S,j,\Gamma,\tilde{u})$ where $(S,j)$ is a
closed Riemann surface, $\Gamma\subset S$ a finite ordered subset,
and $\tilde{u}:S\setminus\Gamma\rightarrow {\mathbb R}\times M$ a
proper map satisfying the partial differential equation
$$
T\tilde{u}\circ j =\tilde{J}\circ T\tilde{u}
$$
and the energy condition
$$
\int_{S\setminus\Gamma} u^{\ast}d\lambda <\infty.
$$
If $\tilde{u}$ solves the differential equation so does
$\tilde{u}_c$ defined by $\tilde{u}=(a+c,u)$, where $c$ is a real
number. We call two such tuples equivalent, say $\alpha \equiv
\alpha'$ provided there exists a biholomorphic map
$\phi:(S,j)\rightarrow (S',j')$ mapping $\Gamma$ onto $\Gamma'$
preserving the ordering so that for a suitable $c$
$$
\tilde{u}'\circ\phi =\tilde{u}_c.
$$
If $S$ is disconnected we only will allow a constant $c$ and not a
locally constant function $c$ in the definition of equivalence. As
it will turn out the invariants of contact manifolds are obtained
by counting equivalence classes of objects (as above) where the
underlying Riemann surface is connected. However, the fact that
such counting leads to invariants independent of the choices
involved also requires counting of certain classes of
non-connected Riemann surfaces. We would like to point out that
even for defining the invariants it is nevertheless important that
non-connected surfaces are being considered. The deeper reason for
this is the following fact. In general it is not possible (for
intrinsic reasons) to pick geometric data in such a way that the
occurring Fredholm operators are transversal to the zero-section.
Nevertheless abstract (multi-valued) perturbations  can be picked
in such a way that transversality can be achieved. In order that
these counts only depend on the connected elements of the moduli
space the perturbations have to be picked in such a way that they
respect the non-connectedness. Observe that due to the fact that
all Fredholm problems are strongly interrelated  there are
infinitely many conditions to be satisfied by the perturbation to
achieve a perturbed problem featuring the same structure. In fact,
due to the relevant compactification which has a certain level
structure, there are connected objects in the compactification
which contain disconnected levels. The perturbation restricted to
such disconnected levels is essentially the same than the
perturbation on the curve in this specific level viewed as part of
a different moduli space.

\subsubsection{Behavior At Punctures}
From now on we assume that all periodic orbits for the Reeb vector
field $X$ are non-degenerate.  A periodic orbit for us will be a
pair $(P,k)$, where $P$ is a submanifold of $M$ diffeomorphic to
$S^1$ and tangent to the Reeb vector field $X$. Moreover $k$ is a
positive integer called the covering number. The following is
known and follows from standard results for Hamiltonian systems:
\begin{proposition}
Given a contact form $\lambda$ there is a Baire set $\Theta\subset
C^{\infty}(M,(0,\infty))$ so that for every $f\in \Theta$ the
contact form $f\lambda$ is non-degenerate.
\end{proposition}

Consider a non-degenerate contact form $\lambda$  on the closed
$M^{2n-1}$. A periodic orbit for the associated Reeb vector field is
a pair $(P,k)$ where $P\subset M$ is diffeomorphic to a circle and
tangent to $X$ and $k\geq 1$ is an integer called the covering
number. We will write $P$ instead of $(P,1)$. The non-degeneracy
means that the linearized Poincare return maps and their iterates do
not have $1$ in the spectrum. The k-fold iterated linearized
Poincare map of $P$ is, of course, the linearized Poincare map for
$(P,k)$ if we take the same local section. Call an orbit $P$
troublesome if $(-1,0)$ contains an odd number of eigenvalues. We
call an orbit $(P,k)$ even if the determinant of $Id-A$ where $A$ is
the linearized Poincare map is positive. Otherwise we call it odd.
Now denote by ${\mathcal P}$ the collection of all period orbits
$(P,k)$ so that $k$ is odd if $P$ is troublesome. We denote the
collection of all periodic orbits by ${\mathcal P}_{all}$. Let us
denote by $\gamma$ the elements of ${\mathcal P}$.

In order to study the behavior of a finite energy map near a
puncture $\gamma\in \Gamma$ fix holomorphic polar coordinates
$\sigma:{\mathbb R}^+\times S^1\rightarrow \dot{S}$ around $\gamma$.
Namely take a disk-like closed neighborhood ${\mathcal D}$ of
$\gamma$ containing no other puncture and having a smooth boundary.
Let $D$ be the closed unit disk in ${\mathbb C}$. Take a
biholomorphic map $h:D\rightarrow{\mathcal D}$ mapping $0$ to
$\gamma$ and define $\sigma$ by
$$
\sigma(s,t)=h(e^{-2\pi(s+it)}).
$$
With $\tilde{u}$ being the finite energy map define
$\tilde{v}=\tilde{u}\circ\sigma$. Then $\tilde{v}$ satisfies
\begin{eqnarray*}
&\tilde{v}:{\mathbb R}\times S^1\rightarrow {\mathbb R}\times M&\\
&\tilde{v}_s+\tilde{J}\tilde{v}_t=0&\\
&E(\tilde{v})<\infty.&
\end{eqnarray*}
Write $\tilde{v}=(b,v)$.

The following result is well-known, see
\cite{Hofer-Weinstein-conj},\cite{HWZ-asymptotics,HWZ-asymptotics-degeneracies}:
\begin{theorem}
Let $\lambda$ be a non-degenerate contact form on the closed
manifold $M$ and assume that $J$ is an admissible complex
multiplication. With $\tilde{v}$ as described above the limit
$$
T:=\hbox{lim}_{s\rightarrow\infty}\ \int_{S^1}\ v(s)^{\ast}\lambda
$$
exists. If $T=0$ the puncture $\gamma$ for $\tilde{u}$ is smoothly
removable. If $T\neq 0$ the number $|T|$ is the period of  a
periodic orbit  of the Reeb vector field $X$ associated to
$\lambda$. Further there exists a constant $d$ and a
$|T|$-periodic orbit $x$ so that
\begin{eqnarray*}
& b(s,t)- Ts -d\rightarrow 0\ \hbox{as}\ s\rightarrow\infty\ \in C^{\infty}(S^1,{\mathbb R})&\\
& v(s,t)\rightarrow x(Tt)\ \in C^{\infty}(S^1,M)\ \hbox{as}\
s\rightarrow\infty.&
\end{eqnarray*}
\end{theorem}
Even more can be said about the convergence near the puncture and
we give more details soon. According to the cases $T<0, T=0, T>0$
we will distinguish negative, removable and positive punctures.

In order to study the behavior near a puncture in more detail we
need special coordinates. In the lemma below we denote by
$\lambda_0$ the standard contact form
\begin{equation*}
  \lambda_0 = d \vartheta +  \sum_{i=1}^{n-1}x_idy_i
\end{equation*}
on $S^1\times {\mathbb R}^{2(n-1)}$ with coordinates $(\vartheta,
x_1,x_2,\dotsc ,x_n,y_1,y_2,\dotsc ,y_{n-1})$.
\begin{lemma}\label{l2}
Let $(M, \lambda)$ be a (2n-1)-dimensional manifold equipped with
a contact form, and let $x(t)$ be a T-periodic solution of the
corresponding Reeb vector field $\dot{x} = X_{\lambda} (x)$ on
$M$. Let $\tau$ be the minimal period such that $T=k\tau$ for some
positive integer k. Then there is an open neighborhood $U \subset
S^1 \times {\mathbb R}^{2(n-1)}$ of $S^1 \times \{0\}$ and an open
neighborhood $V \subset M$ of $P = \{x(t) \mid t \in {\mathbb
R}\}$ and a diffeomorphism $\varphi \colon U \to V$ mapping $S^1
\times \{0\}$ onto $P$ such that
\begin{equation} \label{e2.9}
  \varphi^{\ast} \lambda = f \cdot \lambda_0,
\end{equation}
with a positive smooth function $f \colon U \to {\mathbb R}$
satisfying
\begin{equation} \label{e2.10}
  f\equiv \tau \quad \text{and} \quad df\equiv 0
\end{equation}
on $S^1\times \{0\}$.
\end{lemma}

For large $s$ we can write $\tilde{v}$ in local coordinates around
the limiting periodic orbit as follows (assuming that $T\neq 0$).
$$
\tilde{v}(s,t)= (b(s,t),\vartheta(s,t),z(s,t))
$$
where $z=(x,y)$. Here $s\geq s_0$ and $t\in {\mathbb R}$. The
function $\vartheta(s,t)$ satisfies
$\vartheta(s,t+1)=k+\vartheta(s,t)$. The functions $b$ and $z$ are
$1$-periodic in $t$.

The main result concerning the asymptotic behavior is the following,
\cite{HWZ-asymptotics,HWZ-asymptotics-degeneracies}.
\begin{theorem}
There exists a constant $d>0$ and constants $b_0$, $\vartheta_0$
so that for every multi-index $\alpha$ there is a constant
$C_{\alpha}$ so that
\begin{eqnarray}
& |\partial^{\alpha}[b(s,t)-b_0-Ts]|\leq C_{\alpha} e^{-ds}&\\
& |\partial^{\alpha}[\vartheta(s,t)-\vartheta_0-kt]|\leq C_{\alpha} e^{-ds}&\nonumber\\
&|z(s,t)|\leq C_{\alpha} e^{-ds}&\nonumber
\end{eqnarray}
\end{theorem}

This result is crucial for the functional analytic set-up, since it
tells us which function spaces to take. The constant $d$ is related
to  the spectral properties of some self-adjoint operator $A$
associated to the limiting periodic orbit. In fact $d$ should be
smaller than the distance of the smallest positive or largest
negative eigenvalue to $0$.

\subsection{The Moduli Spaces}
Now we can introduce the moduli spaces we are interested in

\subsubsection{Height-k-Curves} In this subsection we introduce the relevant
compactifications for the moduli spaces of pseudoholomorphic curves
in symplectized contact manifolds. There are obvious extensions to
symplectic cobordisms which are important. We refer the reader to
\cite{BEHWZ,EGH} for complete detail.

The following type of level-k-curves are important and were
introduced as objects occurring in the compactification of moduli
spaces. Let us begin with level-1-curves. Consider tuples
$(S,j,\tilde{u},{\Gamma},M,D)$ where $(S,j)$ is a closed Riemann
surface and $D$ a set of nodal pairs. The sets $\Gamma$ and $M$
are mutually disjoint from the nodal points. Here $M$ is an
ordered set of marked points and ${\Gamma}$ is an ordered set of
so-called punctures. Further
$\tilde{u}:S\setminus\Gamma\rightarrow {\mathbb R}\times M$ is a
$\tilde{J}$-holomorphic proper map and $\tilde{u}(x)=\tilde{u}(y)$
for every nodal pair.

\begin{figure}[htbp]
\mbox{}\\[2ex]
\centerline{\relabelbox
\epsfxsize 3.7truein \epsfbox{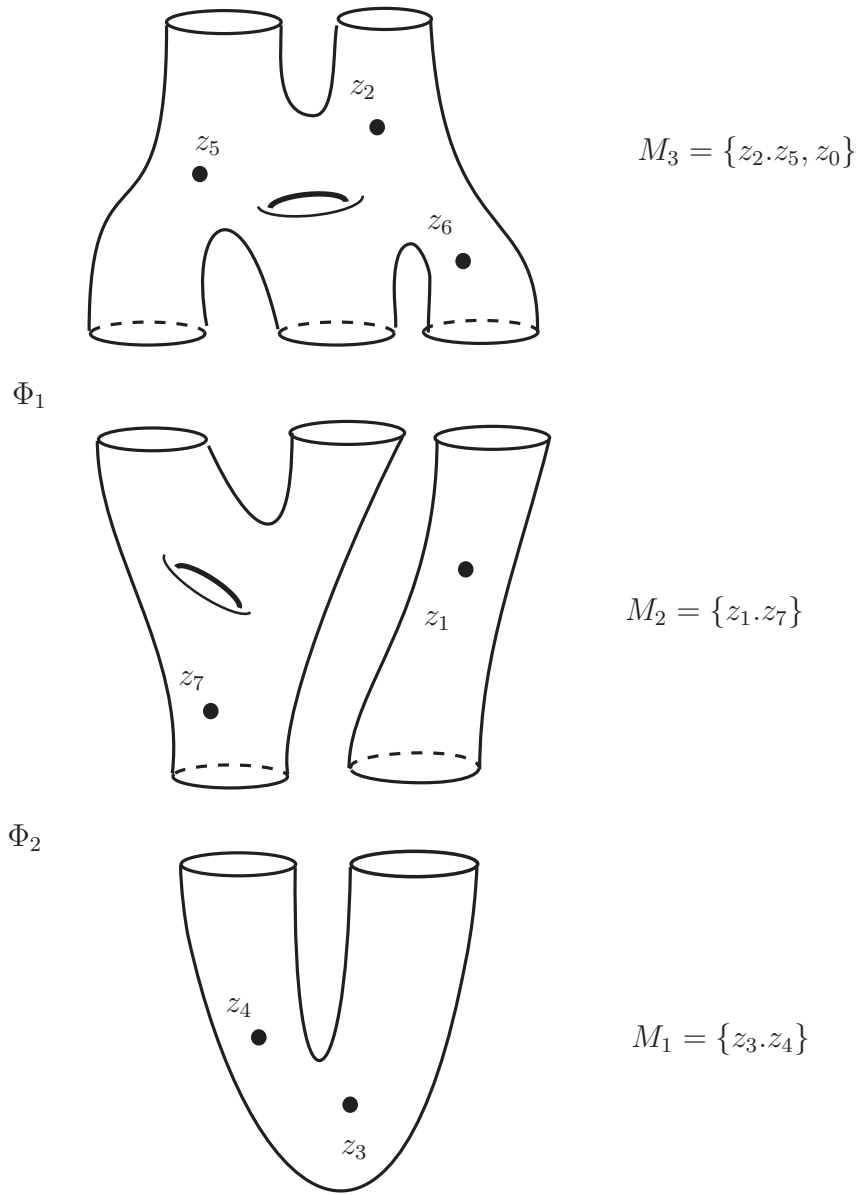}
\relabel {a}{$z_1$}
\relabel {b}{$z_2$}
\relabel {c}{$z_3$}
\relabel {d}{$z_4$}
\relabel {e}{$z_5$}
\relabel {f}{$z_6$}
\relabel {g}{$z_7$}
\relabel {h}{$\Phi_1$}
\relabel {k}{$\Phi_2$}
\relabel {l}{$M_1=\{z_3. z_4\}$}
\relabel {m}{$M_2=\{z_1. z_7\}$}
\relabel {n}{$M_3=\{z_2. z_5, z_0\}$}
\endrelabelbox}
\caption{Holomorphic building of height three with an ordered set
of marked points}\label{Figure5}
\mbox{}\\[1ex]
\end{figure}

 We call two such tuples equivalent if
there exists a biholomorphic map $\phi:(S,j,{\Gamma},M,D)\rightarrow
(S',j',{\Gamma}',M',D')$ and a constant $c$ so that
$\tilde{u}'\circ\phi=\tilde{u}_c$.

\begin{figure}[htbp]
\mbox{}\\[2ex]
\centerline{\relabelbox
\epsfxsize 2.4truein \epsfbox{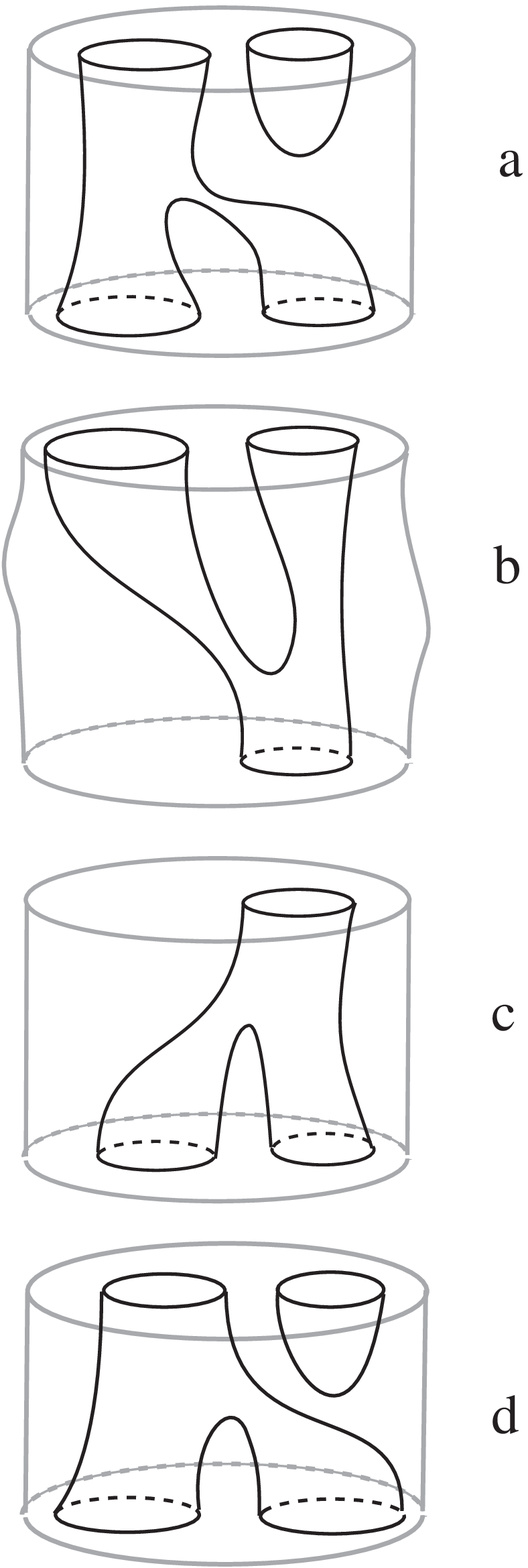}
\relabel {a}{$\R\times V_{+}$}
\relabel {b}{$W$}
\relabel {c}{$\R\times V_{-}$}
\relabel {d}{$\R\times V_{-}$}
\endrelabelbox}
\caption{Holomorphic building of height  $2\vert 1\vert1$}\label{Fig6}
\mbox{}\\[1ex]
\end{figure}

After having introduced level-$1$-curves we define the somewhat more
complicated level-$k$-curves. We are given the following data. A
finite sequence of closed Riemann surface $(S_\ell,j_\ell)$ with
$\ell=1,...,k$. For every $S_{\ell}$ a set of nodal pairs $D_{\ell}$
on $S_{\ell}$ for every $\ell=1,...,k-1$ a set $\hat{D}_{\ell}$ of
ordered pairs $(\hat{x},\hat{y})$ where $\hat{x}$ is a positive
decorated\footnote{A decorated puncture consists of a point $z\in S$
together with an oriented real half line in $T_zS$. The decoration
allows to take a special class of holomorphic polar coordinates
compatible with the asymptotic direction. In particular, the results
about asymptotic convergence imply that there is in an obvious way a
special point on the asymptotic periodic orbit associated to the
asymptotic marker. } puncture on $S_{\ell}$ and $\hat{y}$ a negative
decorated puncture on $S_{\ell +1}$. Moreover an ordered set of
punctures ${\Gamma}$ which is contained in $S_1\cup S_k$ and an
ordered set of marked points $M$ which lie on the union of the
$S_\ell$. The set of $\Gamma$-points on $S_k$ are positive and those
on $S_1$ negative punctures. Moreover we are given maps
$\tilde{u}_{\ell}$ on the $S_{\ell}$. All this data satisfies the
following conditions. At a nodal pair $\{x,y\}\in D_\ell$ we have
$\tilde{u}_\ell(x)=\tilde{u}_\ell(y)$. At the decorated nodes
$(\hat{x},\hat{y})\in \hat{D}_{\ell}$ we require that $\hat{x}$ has
a positive asymptotic limit which is the negative asymptotic limit
at $\hat{y}$ and in addition the points on the periodic orbit
associated to the asymptotic markers coincide. The equivalence of
two such objects is defined via biholomorphic maps between the
corresponding levels preserving all data\footnote{Here we allow the
common rotation of the asymptotic markers of a pair
$(\hat{x},\hat{y})$ beforehand.} and establishing a correspondence
between the maps $\tilde{u}$ where we are allowed a different
${\mathbb R}$-shift on every level. Here is an important definition:
\begin{definition}
We say that a level-k-curve is stable provide for every component
$C$ of a level either $2g +r\geq 0$, where $g$ is the genus of $C$
and $r$ the number of special points on $C$. If that is not the
case we require that the difference of the sum of the periods of
the positive punctures minus the sum of the periods of the
negative punctures is positive.
\end{definition}
Figure \ref{Figure5} shows a level-3-curve in ${\mathbb R}\times V$.
We can also study $\tilde{J}$-holomorphic curves in a symplectic
cobordism with contact-type boundary. Here one adds to the convex
boundary $V^+$ the half-cylinder $[0,\infty)\times V^+$ and to the
concave one $(-\infty,0]\times V^-$. In these necks the analysis is
quite similar to the one described for ${\mathbb R}\times M$. The
analysis in the symplectic cobordism is like Gromov's analysis, see
\cite{Gromov1985,Hummel}. The necessary compactification then
consists of so-called $(k^-|1|k^+)$-curves. An example is depicted
in the Figure \ref{Fig6}. It is almost apparent why this problem
should have something like a module structure over the ${\mathbb
R}\times V^\pm$-problems.

\subsubsection{The polyfold Fredholm set-up}

This compactification of level-1-curves via higher level curves can
be viewed as the zero-set of a polyfold Fredholm problem with
operations. In order to construct a polyfold set-up, we essentially
take the same objects but do not require the maps to be
$\tilde{J}$-holomorphic. The level-structure will incorporate the
differentiability of the maps $\tilde{u}$ in some Sobolev class and
exponential decay properties in their convergence to periodic orbits
near the punctures and some similar properties near the nodes. We
will not describe this in more detail her, but refer the reader to
\cite{EH,HWZ1,HWZ2,HWZ3}. Let us nevertheless describe the indexing
of the operation.

\subsubsection{The operation}

Assuming that we have put the SFT problem into a polyfold set-up
with bundle $Y\rightarrow X$ and a Fredholm section $f$, which , of
course, is the nonlinear Cauchy-Riemann operator, we describe in
this subsection the operation\footnote{There are many variations
about what follows below. Our conventions are somewhat different
from those in \cite{EGH}.}

For every periodic orbit $\gamma$ introduce two symbols
$p_{\gamma}$ and $q_{\gamma}$. These symbols have a grading via a
suitably normalized Conley-Zehnder index (or a mod $2$ reduction
thereof). This index plays the role of a Morse-index. We introduce
a calculus of these symbols by allowing two "even" or one "even"
and one "odd" symbol to commute. Two "odd" symbols anti-commute.
Then we introduce an additional symbol $\hbar$ which is even. We
add the relation
$$
[p_{\gamma},q_{\gamma}] = \kappa_{\gamma}\hbar,
$$
where $\kappa_{\gamma}$ is the covering number of $\gamma$. The
commutator is, of course, a super-commutator. We allow now finite
formal products of powers of these symbols of the form (and in the
order as written)
$$
\tau:=\hbar^{g-1}q_{\gamma_1}^{k_1}...
q_{\gamma_\ell}^{k_{\ell}}p_{\gamma_{\ell+1}}^{k_{\ell+1}}...
p_{\gamma_{\ell+m}}^{k_{\ell+m}}.
$$
The indexing then assigns to $\tau$ the component (not necessarily
connected) of elements in the polyfold which have genus $g$,
asymptotic negative limits (with multiplicities)
$\gamma_1$,...,$\gamma_{k_\ell}$ and asymptotic positive limits
$\gamma_{\ell}$,..,$\gamma_{\ell+m}$. We say this symbol sequence
has standard form. We say two symbols in standard form are
equivalent if by permutation within the $q-$ and $p$-part using the
commutation rules they can be brought into the same form. Denote for
such a symbol sequence $\sigma$ by $[\sigma]$ the equivalence class.
Given $[\sigma]$ and $[\tau]$ define $[\sigma][\tau]$ by
$[\sigma\tau]$. The latter symbol is not in standard form, but using
all the rules is a formal finite sum
$$
[\sigma][\tau]=\sum \lambda_{[\beta]} [\beta],
$$
where the $\beta$ are in standard form and the occurring classes
are different. Further $\lambda_{[\beta]}$ are integers $\neq 0$.
Then define a degeneration structure $(S,R)$ as follows. The set
$S$ consists of all classes $[\sigma]$ in standard form using only
symbols associated to "non-troublesome" periodic orbits, and the
relators are triple $([\sigma],[\tau];[\beta])$, where $[\beta]$
occurs with a nontrivial coefficient in the formal sum above.

The operation then assigns to $[\beta]$ level-k-curves of arithmetic
genus $g$ ($g-1$ is the exponent of $\hbar$) so that the top
punctures are $(+)$-asymptotic to the $\gamma$ occurring in
$p_{\gamma}$ and the bottom punctures are $(-)$-asymptotic to the
$\gamma$ occurring in a $q$-symbol. Also trouble-some orbits will
play a role in the theory. They occur in certain situations which
one might call "geometric wall-crossing".

The homological data then produced from the counting of solutions in
the moduli spaces has a very rich representation theory and we refer
the reader to \cite{EGH}.
\subsection{Comments}
There are in general quite a number of different way to turn moduli
problem into a Fredholm problem with operation. For example in SFT
one might fix asymptotic markers to the punctures and tag the
periodic orbits with generically chosen points, or one might not
make such choices at all. In certain situations this might lead to
moduli spaces which in the first case have M-polyfold descriptions
which are easier, or, in the second case only polyfold descriptions.
In some sense the first description is a covering of the second.
Depending on the situation one might prefer one description over the
other. In \cite{HWZ2} we will describe a certain number of different
approaches to the same problem.

\section{Outlook and Thoughts}
In the cases of Floer-Theory, Gromov-Witten Theory, and SFT there
are many benefits of the polyfold theory. It gives a clean and easy
language to describe and handle these problems. On the abstract
level it offers the benefits of the usual "Fredholm package", i.e.
transversality and perturbation theory. Of course, to bring the
problems into such a framework is usually quite technical (One
shouldn't forget that the problems are for good reason considered
very hard problems, i.e. one should not expect a free ride.) It is
worthwhile to note that the known procedures of bringing the
concrete problems into the polyfold set-up indicates standardized
features. For example it is quite feasible that the understanding of
a wider range of applications would allow us to formulate a certain
number of useful results in a "Scale-Analysis" which would simplify
the transition from a concrete problem to a polyfold description.

It seems to be plausible that the theory described here, or suitable
generalizations, should be applicable to quite a number of nonlinear
problems. In fact a quick look at current research activities shows
that problems with a lack of compactness are very prominent. There
are, of course, problems like Yang-Mills or Seiberg-Witten-Floer
Homology which perhaps (under presumably mild generalizations) could
be put into such a framework. It is likely that the analytical
set-up for proving the Atiyah-Floer conjecture, which requires
ultimately a homotopy from a Yang-Mills to a symplectic Lagrangian
intersection problem, should be possible within our new framework.

Other problems of interest might be Ginzburg-Landau type problems,
or elliptic problems with limiting Sobolev exponents. Here one
should try to derive a good compactification of the problem. Even if
in a physical context only particular solutions might be of interest
it might still be the right point of view to consider a compactified
solution space which carries invariants which cannot be destroyed
and then to show later that for topological reasons there have to be
solutions of physical interest as well\footnote{One should just
recall how useful and fruitful the notion of a weak solution for an
pde has been.}. It seems that currently this type of idea has been
implemented quite successfully in an interesting array of problems
with a geometric background, where the geometry quite often
"dictates" a suitable compactification, but to a lesser degree in
other problems.

Another interesting direction could be concerned with bubbling-off
in a context of geometric evolution problems. It would be
interesting to know if one can describe such phenomena in a polyfold
context. Is there for example a theory of evolution equations in
polyfolds and, importantly, are there good applications which would
show the benefit of such an approach.

Then there are quite a number of "immediate spin-off ideas". For
example as mentioned before one could try to develop some algebraic
topology framework for spaces with operations, i.e. spaces with
boundary with corners, where the faces are explained as products of
their components or more generally as fibered products. In short,
spaces with
$$
\partial M = M\circ M.
$$
Then, as our limited experience already shows there should be some
kind of representation theory of the rudimentary algebraic topology
data leading to interesting algebraic objects. For example
Floer-Theory, Contact Homology and more generally SFT show that the
basic algebraic structures associated to a Fredholm problem with
operations can have representations as differential algebras,
(super-) Poisson  algebras and (super-) Weyl algebras. With other
words the area looks interesting enough to give these issues some
further thought.

\end{document}